\pdfoutput=1

\documentclass[11pt,a4paper]{amsart}    
\usepackage{amsfonts,amsmath,latexsym,amssymb} 
\usepackage{amsthm}                

\usepackage[utf8]{inputenc}
\usepackage[T1]{fontenc}
\usepackage{mathtools}
\mathtoolsset{showonlyrefs=true} 
\usepackage{times}  
\usepackage{xspace}
\usepackage[final]{microtype} 
\usepackage{paralist} 
\usepackage{enumitem}
\usepackage[pdfusetitle,colorlinks,bookmarksopen=true,bookmarksnumbered=true,linkcolor=blue,urlcolor=blue,citecolor=blue]{hyperref}
\usepackage{lineno}
\usepackage{pgfplots}
\pgfplotsset{compat=newest}
\usetikzlibrary{patterns}
\usepgfplotslibrary{fillbetween}
\usepackage{tikzscale}
\usepackage{colortbl}
\usepackage{booktabs}
\usepackage{algorithm}
\usepackage[noend]{algpseudocode}
\usepackage[textsize=tiny,shadow]{todonotes}
 \newtheorem{theorem}{Theorem}[section]
\newtheorem{lemma}[theorem]{Lemma}

\newtheorem{prop}[theorem]{Proposition}
\newtheorem{proposition}[theorem]{Proposition}
\newtheorem{corollary}[theorem]{Corollary}
\newtheorem{hypo}[theorem]{Hypothesis}
\newtheorem{lemmadef}[theorem]{Lemma and Definition}

\theoremstyle{definition}
\newtheorem{definition}[theorem]{Definition}
\newtheorem{example}[theorem]{Example}

\newtheorem{remark}[theorem]{Remark}

\newcommand{\todoinline}[1]{\todo[inline,size=\normalsize]{#1}}
\DeclareMathOperator{\dist}{\mathrm{dist}}
\DeclareMathOperator{\Ran}{\mathrm{Ran}}
\DeclareMathOperator{\Ker}{\mathrm{Ker}}
\DeclareMathOperator{\Kern}{\mathrm{Ker}}
\DeclareMathOperator{\Rank}{\mathrm{Rank}}
\DeclareMathOperator{\spec}{\mathrm{spec}}
\DeclareMathOperator*{\slim}{s-lim}
\DeclareMathOperator*{\wlim}{w-lim}
\newcommand{\cA}{\mathcal{A}}
\newcommand{\cB}{\mathcal{B}}
\newcommand{\cG}{\mathcal{G}}
\newcommand{\cL}{\mathcal{L}}
\newcommand{\cF}{\mathcal{F}}
\newcommand{\dd}{\mathrm d}
\newcommand{\ee}{\mathrm e}
\newcommand{\ii}{\mathrm i}
\newcommand{\eps}{\varepsilon}
\newcommand{\limeps}{\lim_{\eps \rightarrow 0^+} }
\newcommand{\abs}[1]{\left| #1 \right|} 
\newcommand{\scalarproduct}[2]{\langle #1, #2 \rangle}
\newcommand{\norm}[1]{\|#1\|}
\newcommand{\normbig}[1]{\left\|#1\right\|}
\newcommand{\sca}[1]{\langle #1 \rangle}
\newcommand{\scabig}[1]{\left< #1 \vphantom{#1} \vphantom{#1} \right>}
\newcommand{\NN}{\mathbb{N}}
\newcommand{\RR}{\mathbb{R}}
\newcommand{\CC}{\mathbb{C}}
\newcommand{\ZZ}{\mathbb{Z}}
\newcommand{\II}{\mathbb{I}}
\newcommand{\JJ}{\mathbb{J}}
\newcommand{\QQ}{\mathbb{Q}}
\newcommand{\Wcal}{\mathcal{W}}
\newcommand{\ess}{\mathrm{ess}}
\newcommand{\EE}{\mathsf{E}}
\newcommand{\one}{\mathbf{1}}
\newcommand{\boldI}{\mathbf{I}}
\newcommand{\Onetwo}{\mathds{1}_2}
\newcommand{\onematrix}{\mathds{1}}

\newcommand{\Lla}{L_{\lambda,\alpha}}
\newcommand{\Llam}{L_{\lambda,\alpha_m}}

\renewcommand{\P}{\ensuremath{\mathcal{P}}}
\newcommand{\fol}[2]{(#1_{#2})_{#2\in\NN}}
\newcommand{\foln}[2]{(#1)_{#2\in\NN}}
\newcommand{\sfol}[2]{\{#1\}_{#2\in\NN}}
\newcommand{\folz}[2]{(#1_{#2})_{#2\in\ZZ}}
\newcommand{\LL}{\ensuremath{\mathcal{L}}}
\newcommand{\plim}{\operatorname{\P-\lim}\limits}
\newcommand{\pli}{\stackrel{\P}{\rightarrow}}
\newcommand{\Ima}{\textnormal{Im\,}}
\DeclareMathOperator{\supp}{supp}
\DeclareMathOperator{\diam}{diam}

\newcommand{\red}[1]{{\color{red} #1}}
\newcommand{\blue}[1]{{\color{blue} #1}}
\newcommand{\gray}[1]{{\color{gray} #1}}

\newcommand{\eg}{\mbox{e.g.}\xspace}
\newcommand{\ie}{\mbox{i.e.}\xspace}

\newcommand{\energy}{E}
\newcommand{\period}{K}
\newcommand{\tr}{\operatorname{tr}}
\newcommand{\trans}{T}
\newcommand{\Lim}{\operatorname{Lim}}
\newcommand{\cfin}{\mathrm{c}_{00}}
\newenvironment{addmargin}[2][\empty]{\par
  \rightskip=#2\relax
  \ifx\empty#1\relax \leftskip=\rightskip
  \else \leftskip=#1\relax
  \fi}{\par}

\def\drawBand{\fill[bandSpecCol]}
\definecolor{pointSpecCol}{rgb}{0.0, 0.5, 0.0} 
\definecolor{bandSpecCol}{rgb}{0.52, 0.52, 0.51} 
 
\hypersetup{
pdfcreator={LaTeX with hyperref (master, 9719a952)},
pdfkeywords={aperiodic operators, Sturmian words, spectra, limit operators, Jacobi operators},
pdfsubject={MSC Primary 65J10, 47B36; Secondary 47N50},
}

\begin{document}
\title[FSM for aperiodic Schr\"odinger operators]{Finite section method for\\ aperiodic Schr\"odinger operators}

\author[Gabel, Gallaun, Gro{\ss}mann, Lindner, Ukena]{Fabian Gabel, Dennis Gallaun, Julian Gro{\ss}mann,\\Marko Lindner, Riko Ukena} 
\address{ 
Hamburg University of Technology,
Institute of Mathematics, 
Am Schwarzen\-berg-Campus 3, 
D-21073 Hamburg,
Germany} 
\email{\href{mailto:fabian.gabel@tuhh.de}{fabian.gabel@tuhh.de}} 
\email{\href{mailto:dennis.gallaun@tuhh.de}{dennis.gallaun@tuhh.de}} 
\email{\href{mailto:julian.grossmann@jp-g.de}{julian.grossmann@jp-g.de}} 
\email{\href{mailto:lindner@tuhh.de}{lindner@tuhh.de}} 
\email{\href{mailto:riko.ukena@tuhh.de}{riko.ukena@tuhh.de}} 

\date{February 20, 2023}

\begin{abstract}
We consider 1D discrete Schrödinger operators with aperiodic potentials given by a Sturmian word, which is a natural generalisation of the Fibonacci Hamiltonian. Via a standard approximation by
periodic potentials, we establish Hausdorff convergence of the corresponding spectra for the Schrödinger operators on the axis as well as for their compressions to the half-axis.

Based on the half-axis results, we study the finite section method, which is another operator approximation, now by compressions to finite but growing intervals, that is often used to solve operator equations approximately. We find that, also for this purpose, the aperiodic case can be studied via its periodic approximants.
Our results on the finite section method of the aperiodic operator are illustrated by confirming a result on the finite sections of the special case of the Fibonacci Hamiltonian.
\end{abstract}

\maketitle

\section{Introduction}\label{sec:intro}
We study the spectrum and the sequence of so-called finite sections of an aperiodic discrete Schrödinger operator via its periodic approximants. With focus on an expository presentation, we show how to recover results of S\"ut\H{o}~\cite{Suto.1987}, Bellissard~\cite{Bellissard.1989}, and Kellendonk~\cite{Kellendonk.2019} on the spectral approximation on the axis as well as on the half-axis by different means. Precisely, we demonstrate an approach via the method of limit operators and the corresponding toolbox~\cite{Lindner.2006,Rabinovich.2004}, large parts of which are not limited to normal, let alone self-adjoint operators -- see e.g.~\cite{Gabel.2022} for extensions of our results and further details.

Knowledge about the half-axis case for the whole subshift allows us to characterise the stability of the sequence of compressions to finite but growing intervals -- the finite sections -- in both the periodic \cite{Gabel.2021b} and the aperiodic case. The consequences of this analysis generalise results on the finite sections of the Fibonacci Hamiltonian \cite{Lindner.2018} and are illustrated in this special case.

Our presentation is in large parts self-contained, including detailed discussions about the situation in $\ell^p$ with $p\in[1,\infty]$, a precise derivation of the periodic approximants and the steps from zero to one and then two boundary points with homogeneous Dirichlet conditions.

\medskip

Discrete Schrödinger operators are used to describe
physical systems on lattices and, therefore, play an important role in theoretical solid-state physics. One mostly finds two important physically motivated cases here. 
On the one hand, 
a periodic potential function describes a periodic system
such as a periodic tight-binding model.
On the other hand, 
disordered solid media correspond to random potentials. A very interesting
mathematical model between these two cases 
is the one that describes so-called \emph{quasicrystals}, which is an ordered structure without any periodicity in it. They were first found by Shechtman in 1982, see~\cite{Shechtman.1984},
who received the \emph{Nobel Prize in Chemistry 2011} for this discovery.

Quasicrystals can be modelled by Schrödinger operators with aperiodic potentials. The most famous example of
an aperiodic structure for such a two-dimensional quasicrystal
is the so-called \emph{Penrose tiling}, cf.~\cite{Steinhardt.1996}.
Here, we will study a standard one-dimensional quasicrystal model on the spaces $\ell^p(\ZZ)$, $p \in [1, \infty]$, where the \emph{Fibonacci Hamiltonian}, cf.~\cite{Kohomoto.1983, Ostlund.1983}, serves as one well-understood example, see, e.g.,~\cite{Damanik.2016,Lindner.2018}.

Throughout this work, we consider the discrete one-dimensional Schrödinger operator $H_{\lambda,\alpha,\theta} \colon \ell^p(\ZZ) \rightarrow \ell^p(\ZZ)$, $p \in [1, \infty]$, with Sturmian potential, namely,
	\begin{equation}\label{eq:SchrodingerOp}
		(H_{\lambda,\alpha,\theta} x)_n =  x_{n+1} +   x_{n-1} + \lambda v_{\alpha,\theta}(n) x_n \, , \quad n \in \ZZ \,,
	\end{equation}
	where
	\begin{equation*}\label{eq:defintion_words}
v_{\alpha,\theta} (n) = \chi_{[1-\alpha,1)}(n\alpha+\theta \bmod 1)
	\end{equation*}
with coupling constant $\lambda\in\RR$,
irrational slope $\alpha\in [0,1]$, and $\theta \in [0,1)$, cf.~\cite{Luck.1986,Senechal.1995}. 
The already mentioned Fibonacci Hamiltonian arises when setting  $\alpha=\frac 12(\sqrt 5-1)$.

The following picture explains the construction process for the potential~$v_{\alpha,\theta}(n)$ in~\eqref{eq:SchrodingerOp}:

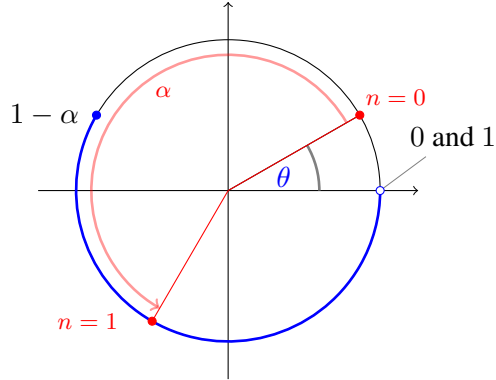
\begin{center}
	\begin{tikzpicture}
		\draw[->] (-2.5,0) -- (2.5,0);
		\draw[->] (0,-2.5) -- (0,2.5);
		\draw (0,0) circle (2cm);
		\draw [-, color=blue, line width =1pt] (150:2) arc (150:360:2cm);
		\draw [-, color=gray, line width =1pt] (0:1.2) arc (0:30:1.2cm);
		\coordinate[pin= 60:{$ 0 $ and $1$}] (A)  at ( 2,0);
		\filldraw[color=blue] (150:2) circle[radius=1.5pt];
		\coordinate[label=left:$1- \alpha ~~$]  (M) at (150:2);
		\coordinate[label=right: {\footnotesize{$\textcolor{red}{\alpha} \quad$}}]  (Y) at (130:1.7);
		\coordinate[] (O)  at ( 0,0);
		\draw[color=blue, fill=white] (0:2) circle[radius=1.5pt];
		\coordinate[label=right:$  \blue{\theta}$] (Tet)  at (0.5,0.18);
		\coordinate[]  (T) at (30:2);
		\draw    (O) -- (T);
		\filldraw[color=red] (30:2) circle[radius=1.5pt];
		\coordinate[label=above: {\footnotesize{$\qquad \quad \textcolor{red}{n = 0}$}}]  (M) at (30:2);
		\draw[color=red]    (O) -- (T);
		\draw [->, color=red, opacity=0.4, line width =1pt] (30:1.8) arc (30:240:1.8cm);
		\filldraw[color=red] (240:2) circle[radius=1.5pt];
		\coordinate[label=left: {\footnotesize{$\textcolor{red}{n = 1} \quad$}}]  (M) at (240:2);
		\draw[color=red]    (O) -- (M);
	\end{tikzpicture}
\end{center}
Since $\alpha$ is irrational, this potential is not periodic.
In the same way, we also consider a similar operator,
	\begin{equation}\label{eq:SchrodingerOpTilde}
            (\widetilde{H}_{\lambda,\alpha,\theta} x)_n =  x_{n+1} +   x_{n-1} + \lambda \widetilde{v}_{\alpha,\theta}(n) x_n \, , \quad n \in \ZZ \,,
	\end{equation}
	where
		\begin{equation*}
\widetilde{v}_{\alpha,\theta} (n) = \chi_{(1-\alpha,1] \cup \{0\} }(n\alpha+\theta \bmod 1)
\,.
	\end{equation*}
Note that, compared to $v_{\alpha,\theta}$, only the interval boundaries have changed in terms of switching from $[\cdot,\cdot)$ to $(\cdot,\cdot]$. Since we identify $0$ and $1$ here, we will always omit mentioning the union with $\{0\}$ from now on.

We will study the following problem in the numerical analysis of aperiodic Schrödinger operators:
Each discrete Schrödinger operator acts on the space~$\ell^p(\ZZ)$ via multiplication with a two-sided infinite tridiagonal matrix~$A$.
If~$A$ is invertible, in order to find the unique solution of the system~$Ax = b$, one often uses a truncation technique which replaces the original (infinite-dimensional) system by a sequence of finite linear systems: the \emph{finite section method (FSM)}. 
In a nutshell, with respect to the standard basis, the operator $A$ and the right-hand side $b$ correspond to a two-sided infinite matrix $(a_{ij})_{i,j\in\ZZ}$ and vector $(\beta_i)_{i\in\ZZ}$, respectively. The FSM replaces $Ax=b$ by $A_nx_n=b_n$, where $A_n=(a_{ij})_{i,j=-n}^n$ and $b_n=(\beta_i)_{i=-n}^n$.
As natural as this method seems at first glance, 
are the projected/truncated equations uniquely solvable (assuming $Ax=b$ is so) and do the solutions
actually approximate the solution $x$ as $n\to\infty$?
If the answers to these questions are positive, we call the FSM for the Schrödinger operator \emph{applicable}.

The applicability of the FSM is strongly related to Fredholm properties of $A$ and its one-sided infinite restrictions $A_+$ and $A_-$.
In numerous works over the last two decades, see, e.g.,~\cite{Bellissard.1989,Damanik.2015,DamanikLenz.1999,Suto.1987}, these operators have been studied and interesting results about their spectra have been found.
In this article, we will use the special tool of \emph{limit operators}, see e.g.,~\cite{Lindner.2006,Rabinovich.2004}, in order to explain and reinterpret some of these well-known results in a framework tailored for our applicability analysis.

\medskip

The following classical approximation of aperiodic Schrödinger operators by periodic operators will form the basis of our approach to the spectrum of $H_{\lambda,\alpha,\theta}$ but also to the applicability of its FSM: 
Expand the irrational slope $\alpha$ as a continued fraction,
\[
\alpha\ = \ \frac{1}{a_1 + \frac{1}{a_2 + \frac{1}{a_3 + \dots}}}
\qquad\text{with}\qquad
a_1, a_2,\ldots\in\NN.
\]
Truncating this infinite expression after $m$ divisions, yields a rational number $\alpha_m \coloneqq \frac{p_m}{q_m}$ that we will refer to as the $m$-th rational approximant to $\alpha$.
Then define the \emph{periodic}(!) Schrödinger operators~${H}_{\lambda,\alpha_m,\theta}$, again by equation~\eqref{eq:SchrodingerOp}. One can see 
that, as $m\to\infty$, not only goes $\alpha_m\to\alpha$ but also ${H}_{\lambda,\alpha_m,\theta}\to{H}_{\lambda,\alpha,\theta}$ in the strong operator topology for $p=2$. 

Using our techniques from \cite{Gabel.2022}, we first show convergence of the spectra of ${H}_{\lambda,\alpha_m,\theta}$ to that of ${H}_{\lambda,\alpha,\theta}$ and, secondly, we show that the same holds for the spectra of their half-axis compressions. Finally, via a study of the FSM of the periodic operators ${H}_{\lambda,\alpha_m,\theta}$, as prepared in \cite{Gabel.2021b}, we conclude a test procedure for the applicability of the FSM of the aperiodic operator, ${H}_{\lambda,\alpha,\theta}$, and easily confirm, as a special case, the main result of \cite{Lindner.2018} about the FSM of the Fibonacci Hamiltonian.

\medskip

The paper is organised as follows. In Section~\ref{sec:fsm_intro}, we introduce the reader to the finite section method and prove several auxiliary results that will be needed in the sections to follow.
Next, we give an overview about Sturmian potentials in Section~\ref{sec:sturm}
and also go into detail for continued fraction expansions of irrational numbers, thereby, laying the groundwork for the spectral and approximation theory of aperiodic Schrödinger operators to be presented in subsequent sections.
Finally,
in Section~\ref{sec:specHltheta} and~\ref{sec:fsm}, we analyse the spectrum of the aperiodic Schrödinger operator $H_{\lambda,\alpha,\theta}$, characterise applicability of the FSM, and derive a numerical method to verify the applicability.

\section{Approximation Methods for Band Operators}\label{sec:fsm_intro}
In order to motivate the finite section method, let us consider a discrete Schrödinger operator $H \colon \ell^p(\ZZ) \rightarrow \ell^p(\ZZ)$
for a fixed parameter $p \in [1, \infty]$
given by
\begin{equation*}
	(H x)_n = x_{n+1} +   x_{n-1} 	 + v(n) x_n \, ,\qquad n\in\ZZ \,,
\end{equation*}
        where the \emph{potential} $v$ is a real-valued function defined on $\ZZ$.

\subsection{Motivation of the Finite Section Method}
The operator~$H$ can analogously be described by a two-sided infinite tridiagonal matrix~$A = (a_{ij})_{i,j \in \ZZ}$ which is defined by
	\begin{equation}
          \label{eq:twoSidedInfiniteMatrixRep}
		A = 
	      \begin{pmatrix}
              \ddots & \ddots  &        &      \\	
		\ddots  & 	v(0)      & 1   &        &        &            \\
		&	1     & v(1)    &  1   &        &              \\
         &   &  1    & v(2)    & \ddots &                 \\
           &        &  & \ddots & \ddots         
		\end{pmatrix} \,.
	\end{equation}

The so-called \emph{finite section method (FSM) for~$H$ or $A$} considers the sequence of finite submatrices 
\begin{equation}
  \label{eq:submatricesFSM}
A_n = (a_{ij})_{i,j=-n}^{n},\quad n\in\NN \,,
\end{equation}
and asks the following: 
\vskip1ex
\begin{addmargin}[2em]{2em}
  \noindent\emph{Are the matrices $A$ and $A_n$ invertible for all sufficiently large $n$ and are their inverses (after embedding them into a two-sided infinite matrix) strongly convergent to the inverse of $A$?}
\end{addmargin}
\vskip1ex
In the case of a positive answer, we call the FSM (for $A$) \emph{applicable}.
This approximation of $A^{-1}$ can be used for solving equations $Ax=b$ approximately via the solutions of growing finite systems $A_n$.
First rigorous treatments of this natural approximation can be found in~\cite{Baxter.1963, Gohberg.1974}. 
In case that the FSM is not applicable in the above sense, it may still be possible to establish an applicability result by passing to suitable 
modifications $A_n = (a_{ij})_{i,j=l_n}^{r_n}$ with well-chosen sequences $l_n\to-\infty$ and $r_n\to\infty$, see~\cite{Lindner.2010, Rabinovich.2004, Rabinovich.2008}.

Neither the FSM nor the following tools and results in Chapter \ref{sec:fsm_intro} require self-adjointness (unless we explicitly say so).

Operators like the Schrödinger operator $H$ whose infinite matrix representation only exhibits finitely many non-zero diagonals are a well-known subject of investigation regarding the applicability of the FSM.
Operators of this type are summarised in the class of \emph{band operators}, which we introduce next. 
We follow the definitions given in~\cite[Section~1.3.6]{Lindner.2006}. 

For nonempty $U,V\subset\ZZ$, we define $\dist(U,V)\coloneqq \min \{|u-v| : {u\in U, v\in V}  \}$ 
  and the projection $P_U: \ell^p(\ZZ)\rightarrow\ell^p(\ZZ)$ by 
  \begin{equation*}
  (P_U x)_n \coloneqq \begin{cases}
    x_n & \text{ if } n\in U \,,\\ 
    0 & \text{ if } n\notin U \,.
  \end{cases} 
\end{equation*}

\begin{definition}[Band-width and band operator]\label{def:band-width}
  Let $A$ be a bounded operator on $\ell^p(\ZZ)$ with $p \in [1,\infty]$. 
  We say that $A$ has \emph{band-width less than $w\in\NN$} if~$P_U A P_V = 0$ for all $U, V \subset \ZZ$ with $\dist(U, V ) > w$, and in this case,
  $A$ is called a \emph{band operator}.
\end{definition}

Note that the notion \emph{band operator} makes also sense for unbounded operators on $\ell^p(\ZZ)$ but, in this work, all band operators will be bounded.

By definition, a band operator $A$ with band-width less than $w$ is only supported on the $k$-th diagonals with $|k|\le w$.
In the following, we will identify an operator $A$ on $\ell^p(\ZZ)$ with its usual \emph{matrix representation} $(a_{ij})_{i,j \in \ZZ}$ with respect to the canonical basis in $\ell^p(\ZZ)$ as we did before.

\begin{remark}\label{rem:band-width}
Whenever necessary, we consider for some index set $\II \subset \ZZ$ the space $\ell^p(\II)$ as a subspace of $\ell^p(\ZZ)$.
Then, the terms in Definition~\ref{def:band-width} naturally carry over to operators $A$ on $\ell^p(\II)$  by associating $A$ with its restriction $A P_{\II}^{} \colon \ell^p(\ZZ) \to \ell^p(\ZZ)$. 
\end{remark}

Assuming invertibility of $A$ on $\ell^p(\ZZ)$, the applicability of the FSM is equivalent to the uniform boundedness of the inverses $A_n^{-1}$, a concept that is also known as \emph{stability}.
\begin{definition}[Stability]\label{def:stability}
    A sequence of operators $(A_n)_{n\in\NN}$ defined on $\ell^p(\ZZ)$ for $p \in [1,\infty]$ is called \textit{stable} if there exists~$n_0\in\NN$ such that for all $n\ge n_0$ the operators $A_n$ are invertible and their inverses are uniformly bounded.
\end{definition}
The basic result connecting the notions of applicability and stability is known as \emph{Polski's~theorem}, cf.~\cite[Theorem~1.4]{Hagen.2001},
and it says that $(A_n)$ is applicable to $A$ if and only if it is stable and $A$ is invertible.

Now, stability, and hence, applicability, of the FSM sequence $(A_n)$ in~\eqref{eq:submatricesFSM} is closely connected to the following entrywise limits
\begin{equation} \label{eq:ass_matrix}
(a_{i+l_n,j+l_n})_{i,j=0}^\infty\ \xrightarrow{n\rightarrow \infty}\ B_+\quad\textrm{and}\quad (a_{i+r_n,j+r_n})_{i,j=-\infty}^0\ \xrightarrow{n\rightarrow \infty}  C_-
\end{equation}
of one-sided infinite submatrices of $A$, where we 
consider sequences $(l_n)$ and $(r_n)$ with $\lim_{n \to \infty} l_n = -\infty$ and $\lim_{n \to \infty} r_n = \infty$, such that the limits in~\eqref{eq:ass_matrix} exist. 
The following result summarises the aforementioned connections and is a standard outcome regarding the applicability of the FSM for band operators. 

\begin{lemma}[{\cite[Lemma 1.2]{Chandler-WildeLindner.2016}}]\label{lem:FSMappl}
Let $A$ be a band operator on $\ell^p(\ZZ)$ for $p\in[1,\infty]$.  
Then the following are equivalent:
\begin{enumerate}[label=(\roman*)]
  \item\label{it:FSMAppl} The FSM is applicable.
  \item\label{it:seqStable} The sequence $(A_n)$ with $A_n = (a_{ij})_{i,j=-n}^{n}$ is stable.
\item\label{it:AandLimOpsInv} The operator $A$ and the limits $B_+$ and $C_-$ from~\eqref{eq:ass_matrix} are invertible for all suitable sequences $(l_n)$ and $(r_n)$. 
\end{enumerate}
\end{lemma}

\subsection{Invertibility and Limit Operators}\label{sec:invAndLimOps}
We now focus on an operator-theoretical tool in order to ensure condition~\ref{it:AandLimOpsInv} of Lemma~\ref{lem:FSMappl}: the so-called limit operators, cf.~\cite{Lindner.2006,Rabinovich.1998,Rabinovich.2004}.

\begin{definition}[Limit operators and compressions]\label{def:limOps}
  Let $A$ be a band operator on $\ell^p(\mathbb{I})$ for $p\in[1,\infty]$ and $\mathbb{I}\in\{\NN,\ZZ\}$. 
  An operator $B\in\ell^p(\ZZ)$ with matrix representation $(b_{ij})_{i,j \in \ZZ}$ is called \emph{a limit operator of $A$} if there is a sequence  $h=(h_n)_{n\in\NN} \subset \ZZ$ with 
  $\lim\limits_{n\rightarrow\infty} |{h}_n| = \infty$ and
  \[ a_{i+h_n,j+h_n}\xrightarrow{n \rightarrow \infty} b_{ij} \]  
  for all $i,j\in\ZZ$.
  In this case, we also write $A_h  \coloneqq B$ and say that $h$ is the \emph{corresponding sequence to $B$}. 
  We denote the set of all limit operators of $A$ by $\Lim(A)$.
  For~$A_h\in\Lim(A)$, we write~$A_h\in\Lim_+(A)$ or $A_h\in\Lim_-(A)$ if the corresponding sequence $h$ tends to $+\infty$ or $-\infty$, respectively.
  Moreover, if $\II = \ZZ$, we call both operators
  \begin{equation*}
   A_+ \coloneqq (a_{ij})_{i,j = 0}^\infty \quad\text{and}\quad A_- \coloneqq (a_{ij})_{i,j = -\infty}^0
  \end{equation*}
  \emph{(one-sided) compressions of} $A$. Then we have the identities $\Lim_\pm(A) = \Lim(A_\pm)$.
\end{definition} 

In fact, the operators $B_+$ and $C_-$ from~\eqref{eq:ass_matrix} correspond to the one-sided compressions of the limit operators $A_{(l_n)}$ and $A_{(r_n)}$.
Evidently, a limit operator $A_h$ of a band operator $A$ is again a band operator with the same bandwidth.

We now reformulate Lemma~\ref{lem:FSMappl} in the language of limit operators and, therefore, from now on, we will use the index sets $\ZZ_-, \ZZ_+$, which both include zero,
and the index set $\NN$, which excludes zero.

\begin{prop}\label{prop:FSMCondition}
Let $A$ be a band operator on $\ell^p(\ZZ)$ for $p\in[1,\infty]$.
Then, the FSM for~$A$ is applicable if and only if 
\begin{enumerate}[label=(\alph*)]
\item\label{it:AinvLp} $A$ is invertible on $\ell^p(\ZZ)$,
\item for all $B\in\Lim_+(A)$ the compressions $B_-$ are invertible on $\ell^p(\ZZ_-)$, and 
\item for all $C\in\Lim_-(A)$ the compressions $C_+$ are invertible on $\ell^p(\ZZ_+)$.
\end{enumerate}
\end{prop}

The restriction to a particular choice of $p$ in Proposition~\ref{prop:FSMCondition} can be dropped  without loss of generality as the following proposition shows. 
In subsequent parts of this article, this result will allow us to fall back onto the Hilbert space case~$p = 2$ whenever needed.

\begin{proposition}
    \label{prop:kurbatov}
    Let $B$ be a band operator on $\ell^p(\II)$ for some $\II \subset \ZZ$ and $p \in [1,\infty]$.
    If~$B$ is invertible, then $B$ is also invertible as an operator on $\ell^q(\II)$ for all $q \in [1,\infty]$.
\end{proposition}
\begin{proof}
	First, we identify the Banach space $\ell^p(\ZZ)$ with the $p$-direct sum $\ell^p(\JJ) \oplus \ell^p(\II)$
	where $\JJ \coloneqq  \ZZ \setminus \II$
	and the norm is given by 
        \begin{equation*}
          \norm{x \oplus y} \coloneqq  
          \begin{cases}
            \big( \norm{x}^p_{\ell^p(\JJ)} 
            + \norm{y}^p_{\ell^p(\II)})^{\frac{1}{p}} &\quad\text{if } p < \infty\,,\\[0.5em]
            \max \{
            \norm{x}_{\ell^\infty(\JJ)} , 
            \norm{y}_{\ell^\infty(\II)}  \}  &\quad\text{if } p  = \infty \,.
          \end{cases}
        \end{equation*}
    Now, assume that~$B$ is invertible on $\ell^p(\II)$.
    We extend the operator~$B$ to an invertible operator~$A$ on~$\ell^p(\ZZ)$ via
    \begin{equation*}
        A \coloneqq  
        \begin{pmatrix}
            \one_{\JJ} & 0 \\
            0 & B
        \end{pmatrix}
    \end{equation*}
    with respect to the direct decomposition $\ell^p(\ZZ)  = \ell^p(\JJ) \oplus \ell^p(\II)$.
    Here, $\one_{\JJ}$ denotes the identity on $\ell^p(\JJ)$.
    Note that, as a band operator, $A$ is an element of the so-called \emph{Wiener algebra} $\Wcal$, see~\cite[Definition~1.43]{Lindner.2006}.
    Due to~\cite[Corollary~2.5.4]{Rabinovich.2004} on the inverse closedness of $\Wcal$, 
    the inverse $A^{-1}$ is in $\Wcal$ again and, hence, acts boundedly on every space $\ell^q(\ZZ)$, $q \in [1,\infty]$, so that
    the operator~$A$ is invertible on all spaces~$\ell^q(\ZZ)$, $q \in [1,\infty]$.
    Moreover, we have
    \begin{equation*}
        A^{-1} 
        = 
        \begin{pmatrix}
            \one_{\JJ} & 0 \\
            0 & B^{-1}
        \end{pmatrix} \,.
    \end{equation*}
    Consequently, $B^{-1}$ is a bounded operator on $\ell^q(\II)$ for all $q \in [1,\infty]$.
\end{proof}

Recall that a bounded linear operator $A \colon X \rightarrow Y$
between Banach spaces $X$ and $Y$
is called a \emph{Fredholm operator} if
the kernel $\Ker(T)$
and the cokernel $Y/\Ran(T)$ are finite-dimensional.
Furthermore, we define the \emph{essential spectrum of $A$} by
\begin{equation*}
  \sigma_\ess(A) \coloneqq  \{ E \in \CC  : A - E  \text{ is not a Fredholm operator} \}\,.
\end{equation*}

The next lemma establishes the equivalence of the Fredholmness of a band operator $A$, the invertibility of its limit operators, see~\cite{Lindner.2014,Rabinovich.1998}, or, alternatively, their injectivity on the space $\ell^\infty(\mathbb{Z})$, see~\cite{Wilde-Lindner.2008,Chandler-WildeLindner.2011}.
The lemma previously appeared in this form in~\cite{Lindner.2018}.

\begin{lemma}[{\cite[Lemma 2.6]{Lindner.2018}}]\label{lem:fundamental}
    Let $\II \in \{\ZZ, \ZZ_+, \ZZ_-, \NN\}$ and $p \in [1,\infty]$.
    For a band operator~$A$, the following are equivalent:
    \begin{enumerate}[label=(\roman*)]
      \item\label{it:AFredholm} $A$ is a Fredholm operator on $\ell^p(\II)$.
      \item\label{it:allLimOpsInvertible} All limit operators of $A$ are invertible on $\ell^p(\ZZ)$.
      \item\label{it:allLimOpsInjective} All limit operators of $A$ are injective on $\ell^\infty(\ZZ)$.
    \end{enumerate}
\end{lemma}

\begin{remark}
Actually, it is Lemma \ref{lem:fundamental}, whose first versions go back to \cite{Rabinovich.1998} and somehow finalise in \cite{Wilde-Lindner.2008,Lindner.2014}, where limit operators naturally enter the scene, and they only reappear in theorems on applicability and stability of operator sequences because the latter is equivalent to Fredholmness of an associated operator.
\end{remark}

Lemma~\ref{lem:fundamental} is particularly useful when dealing with a band operator $A$ that additionally has the property~$A \in \Lim(A)$. 
All subsequent examples of Schrödinger operators will have this property, which is sometimes referred to as~\emph{self-similarity}, cf.~\cite{Beckus.2018,Wilde-Lindner.2008}.
The combination of Proposition~\ref{prop:kurbatov} with Lemma~\ref{lem:fundamental} leads to the following results that allow us to translate the invertibility problem of an operator into an injectivity problem of its limit operators.

\begin{corollary}\label{cor:fundamental}
    Let $A$ be a band operator with $A \in \Lim(A)$.
    Then the following are equivalent:
    \begin{enumerate}[label=(\roman*)]
        \item\label{it:LimInj} All limit operators of $A$ are injective on $\ell^\infty(\ZZ)$.
        \item\label{it:AInvAll} $A$ is invertible on~$\ell^p(\ZZ)$ for all  $p \in [1,\infty]$.
        \item\label{it:AInvSome} $A$ is invertible on~$\ell^p(\ZZ)$ for some $p \in [1,\infty]$.
    \end{enumerate}
\end{corollary}

\begin{proof}
  The equivalence \ref{it:AInvAll}$\Leftrightarrow$\ref{it:AInvSome}
  is immediately given by Proposition~\ref{prop:kurbatov}.
  Since $A \in \Lim(A)$, the implication \ref{it:LimInj}$\Rightarrow$\ref{it:AInvSome}
  follows from Lemma~\ref{lem:fundamental}.
  Using that an invertible operator is a Fredholm operator, we also have~\ref{it:AInvSome}$\Rightarrow$\ref{it:LimInj} by Lemma~\ref{lem:fundamental}.
\end{proof}

Note that Corollary~\ref{cor:fundamental} only handles operators that are defined on $\ell^p(\ZZ)$.
However, Lemma~\ref{lem:FSMappl} also relies on the invertibility of one-sided compressions.
The next corollary gives the corresponding result.
For this, we consider only the case
$p = 2$ because, there, we are able to use the \emph{Hilbert space adjoint} $A^*$
for an operator~$A$ on the Hilbert space $\ell^2(\II)$ with $\II \subset \ZZ$.

\begin{corollary}\label{cor:injImpliesInvertibleForLimOps}
    Let $A$ be a self-adjoint invertible band operator on $\ell^2(\ZZ)$ and $B \in \Lim(A)$.
    \begin{enumerate}[label=(\alph*)]
      \item\label{it:B-inj} 
        If the compression $B_-$ is injective on $\ell^\infty(\ZZ_-)$, then $B_-$ is invertible on~$\ell^p(\ZZ_-)$ for all $p \in [1,\infty]$.
      \item\label{it:B+inj} 
        If the compression $B_+$ is injective on $\ell^\infty(\ZZ_+)$, then $B_+$ is invertible on~$\ell^p(\ZZ_+)$ for all $p \in [1,\infty]$.
    \end{enumerate}
\end{corollary}

\begin{proof}
  We only prove~\ref{it:B-inj} since the proof of~\ref{it:B+inj} works completely analogously.
    For the proof of~\ref{it:B-inj}, we note that, due to Proposition~\ref{prop:kurbatov}, it suffices to consider $p = 2$. 
     As $B_-$ is injective on $\ell^\infty(\ZZ_-)$, it is also injective on the subset $\ell^2(\ZZ_-) \subset \ell^\infty(\ZZ_-)$.
    We will show that $B_-$ has a dense and closed range making $B_-$ also surjective. 
    
    \emph{The range of $B_-$ is dense in $\ell^2(\ZZ_-)$:}
    As $A$ is self-adjoint, all limit operators of $A$ are self-adjoint, see~\cite[Proposition~3.4\,e)]{Lindner.2006}.
    In particular, $B$ and its compression $B_-$ are self-adjoint.
    Therefore the adjoint $(B_-)^*$ is also injective on $\ell^2(\ZZ_-)$ which implies that the range of $B_-$ is dense in $\ell^2(\ZZ_-)$.

  \emph{The range of $B_-$ is closed in $\ell^2(\ZZ_-)$:}
    By assumption, $A$ is also a Fredholm operator on~$\ell^p(\ZZ)$.
    By the implication~\ref{it:AFredholm}$\Rightarrow$\ref{it:allLimOpsInvertible} in Lemma~\ref{lem:fundamental}, with $p = 2$, all limit operators of $A$ are invertible on $\ell^2(\ZZ)$. 
    In particular, $B$ is invertible on $\ell^2(\ZZ)$.
    As this implies that $B$ is Fredholm on $\ell^2(\ZZ)$,  Lemma~\ref{lem:fundamental} again gives that all operators in $\Lim(B) \supset \Lim(B_-)$ are invertible on $\ell^2(\ZZ)$.
    As all limit operators of $B_-$ are invertible on $\ell^2(\ZZ)$,
    the implication~\ref{it:allLimOpsInvertible}$\Rightarrow$\ref{it:AFredholm} in Lemma~\ref{lem:fundamental} gives that $B_-$ is Fredholm. 
    In particular, $B_-$ has a closed range. \qedhere
\end{proof}

We close this first part of the section about limit operators with an important result regarding (essential) spectra of band operators and their limit operators, see also~\cite[Corollary 12]{Lindner.2014}.

\begin{proposition}\label{prop:essSpecEqualsSpecnew}
	Let $A$ be a band operator on $\ell^p(\II)$ with $\II \in \{\ZZ, \ZZ_+, \ZZ_-, \NN\}$.
	Then
	\begin{equation*}
		\sigma_\ess(A) = \bigcup_{B \in \Lim(A)} \sigma(B) \, .
	\end{equation*}
	In particular, we have:
		\begin{enumerate}[label=(\alph*)]
			\item If $\sigma(B) = \sigma(B^\prime)$ for all $B, B^\prime \in \Lim(A)$, then $\sigma_\ess(A) = \sigma(B)$.
			\item If $A \in \ell^p(\ZZ)$ and $A \in \Lim(A)$, then $\sigma(A) = \sigma_\ess(A)$.
		\end{enumerate}
\end{proposition}

\begin{proof}
	By definition,
	$\energy \in \sigma_\ess(A)$ if and only if $A-\energy$ is not Fredholm.
	By the equivalence~\ref{it:AFredholm}$\Leftrightarrow$\ref{it:allLimOpsInvertible}
	in Lemma~\ref{lem:fundamental}, this is equivalent to $B-\energy$ not being invertible, i.e.\ $\energy \in \sigma(B)$, for some $B \in \Lim(A)$.
\end{proof}

\subsection{Stability of Approximation Sequences and Invertibility}
Now we analyse the stability of an approximation sequence $(A_n)$ more closely. 
Keeping in mind the introductory example of the FSM for a bounded operator $A \in \ell^p(\ZZ)$, we now consider more general sequences of approximants $(A_n)$ with $A_n x \to A x$ for all $x \in \ell^p(\ZZ)$. We refer to this as \textit{strong convergence} (or convergence in the strong operator topology) and simply write $A_n\to A$.

A convenient tool to study the norm of $A^{-1}$ or $A_n^{-1}$ (which we clearly have to, for stability) is the so-called lower norm. We call it so, by abuse of notation (it's \emph{not} a norm), just like \cite[Definition~2.31]{Lindner.2006}, \cite[Definition~7.1.6]{Rabinovich.2004}, and others do. In Hilbert space, while $\|A\|$ is the largest singular value of $A$, $\nu(A)$ is the smallest.

\begin{definition}[Lower norm]\label{def:lowerNorm}
  Let $\II \subset \ZZ$ and $A$ be a bounded operator on $\ell^p(\II)$. 
  We define the \textit{lower norm of $A$} by
\[
  \nu(A) \coloneqq  \inf \big\{\|Ax\| : \|x\|=1 \big\} \,.
\]
Whenever necessary, we write $\nu^p(A) \coloneqq \nu(A)$ to clarify which norms are used.
\end{definition}

The connection between $\nu(A)$ and $\|A^{-1}\|$ is fairly straightforward:

\begin{lemma}[{\cite[Lemma~2.10]{Lindner.2016}}, {\cite[Lemma~2.35]{Lindner.2006}}]\label{lem:lowerNormInvertible}
Let $p\in[1,\infty]$, $\II\subset \ZZ$
 and $A$ be a bounded operator on $\ell^p(\II)$. Then it holds that
\begin{equation*}
\|A^{-1}\|^{-1}=\min \big\{\nu(A),\nu(A')\big\}.
\end{equation*}
 In particular, $A$ is invertible if and only if $\nu(A)>0$ and $\nu(A')>0$. 
 In that case, $\nu(A)=\nu(A')$, so that $\|A^{-1}\|^{-1} = \nu(A)$.
\end{lemma}
Above, we use the conventions that $\|A^{-1}\|=\infty$ if $A$ is not invertible, that $\frac{1}{\infty} \coloneqq  0$, and that $A'$ is the \emph{dual operator} of $A$. If $p<\infty$, $A'$ acts on $\ell^{p'}(\II)$ with $p'$ referring to the \emph{conjugated Hölder exponent} of $p$, \ie $\frac{1}{p} + \frac{1}{p'} = 1$. Interpreting $A'$ for $p=\infty$ is more difficult. We start with a little technical build-up that we will need at some point anyway:

\begin{definition}[{Support, diameter, finite sequences and null sequences}] \label{def:support}
  For a sequence $x\in \ell^p(\ZZ)$, we define its \emph{support} by $\supp(x) \coloneqq  \{k\in\ZZ : x_k\neq 0\}$ and,
  for any $U\subset \ZZ$, we call $ \diam(U) \coloneqq  \sup_{k,j\in U} |k-j|$ its \emph{diameter}.
Then, for $\II \subset \ZZ$, the space of all finite sequences on $\II$ 
  is denoted by
  \begin{equation*}
    \cfin(\II) \coloneqq \big\{x \in \ell^\infty(\II) : \diam(\supp(x)) < \infty \big\}\,
  \end{equation*}
  and its closure, the set of all null sequences on $\II$, is denoted by
  \begin{equation*}
    c_0(\II) \coloneqq {\mathop{\rm clos}}_{\ell^\infty(\II)} \cfin(\II).
  \end{equation*}
\end{definition}
Now, the dual space of $\ell^\infty(\II)$ contains $\ell^1(\II)$ but also a lot more (that we do not want to study). Its subspace $c_0(\II)$, though, has exactly $\ell^1(\II)$ as its dual space. Now, if $A$ is a band operator and we let it act on $\ell^\infty(\II)$, it leaves the subspace $c_0(\II)$ invariant, so that its restriction $A_0\coloneqq A|_{c_0(\II)}$ acts $c_0(\II)\to c_0(\II)$. We will then interpret $A'$ from Lemma \ref{lem:lowerNormInvertible} as the dual operator $(A_0)':\ell^1(\II)\to\ell^1(\II)$. As it turns out, $((A_0)')'=A$, which is why $(A_0)'$ is also called \emph{predual operator of $A$}, 
cf.\ \cite[Section~1.3.2]{Lindner.2006}.

Maybe it is interesting in its own right (see \cite[Lemma~3.8]{Lindner.2016}) that $A_0$ and $A$ have the same norm, lower norm, spectrum, and that $(A_0)^{-1}=(A^{-1})_0$.

So, thanks to Lemma \ref{lem:lowerNormInvertible}, we are down to the lower norms of $A$ and $A'$. For self-adjoint operators (and $p=2$, or not) this will simplify further.

Recall that a band operator that is bounded  on $\ell^p(\II)$ for some $p\in[1,\infty]$, 
is also bounded on $\ell^p(\II)$ for all $p\in[1,\infty]$. If the two band operators, $A$ and its dual $A'$, agree on one $\ell^p(\II)$ (say, $A=A'$ on $\ell^2(\II)$), they will agree on all $\ell^p(\II)$. 
In this case, we will say that $A$ is \emph{its own dual operator}.

For convenience, we formulate the following immediate result:

\begin{lemma}\label{lem:approxMerged}
Let $\II\subset\ZZ$ and $p\in[1,\infty]$.
Let further $A$ and $(A_n)_{n\in\NN}$ be band operators on $\ell^p(\II)$ that are their own dual operators.
Additionally, let $\nu^q(A_n)\xrightarrow{n \rightarrow \infty }\nu^q(A)$ for all $q\in \{p, p'\}$.
 Then $A$ is invertible if and only if $(A_n)$ is stable. In that case, 
  \begin{equation*}
  \|A^{-1}\| = \lim_{n \rightarrow \infty} \|A_{n}^{-1}\|\,.
  \end{equation*}
\end{lemma}

\begin{proof}
The invertibility of $A$ is equivalent to $\nu^p(A)>0$ and $\nu^{p'}(A')>0$ by Lemma~\ref{lem:lowerNormInvertible}.
By assumption $\nu^p(A_n)\rightarrow\nu^p(A)$ and $\nu^{p'}(A_n')\rightarrow\nu^{p'}(A')$.
Hence, using Lemma~\ref{lem:lowerNormInvertible} again, we get for sufficiently large $n$ that $A_n$ is invertible. The identity $\|A^{-1}\|^{-1}=\nu(A)$ gives that $(A_n)$ is also stable and 
$  \|A^{-1}\| = \lim_{n\rightarrow \infty}\|A_{n}^{-1}\|$
  holds.
  Similarly, the other direction follows.
\end{proof}

For band operators, the following lemma provides a useful tool to analyse the lower norm via the approximate lower norms. 

\begin{lemmadef}[{\cite[Proposition 6]{Lindner.2014}}]\label{lem:lowerNorm} 
  Let $p\in [1,\infty]$ and $\II\subset\ZZ$.  
  Furthermore, let $\varepsilon >0, r>0$ and $w\in\NN$.
  Then there is a constant $D\in\NN$ such that, for all band operators $A$ on $\ell^p(\II)$ with $\|A\| < r$ and band-width less than $w$, we have:
  \begin{equation*}
  \nu(A) \le \nu_D(A) \le \nu(A)+\varepsilon \,,
  \end{equation*}
  where we set for $D\in\NN$
  \[
    \nu_D (A) \coloneqq  \inf\Big\{ \|Ax\| :  \|x\|=1 ,\;  \diam(\supp(x))\le D \Big\} \;
  \]
  and call it the \emph{approximate lower norm} (for diameter $D$) of $A$.
\end{lemmadef}

The following result is an improvement of~\cite[Corollary 2.37]{Lindner.2006} for band operators.

\begin{lemma}\label{lem:stableImpliesInv}
Let $\II\subset \ZZ$ and $p\in[1,\infty]$.
 Further let $A$ and $\fol{A}{n}$ be band operators on $\ell^p(\II)$ with $A_n x \rightarrow Ax$ and $A_n'y \rightarrow A'y$ for all $x,y\in \cfin(\II)$. 
  If $(A_n)$ has a stable subsequence, then $A$ is invertible and
\begin{equation}\label{eq: lim inf norm inverse estimate}
  \|A^{-1}\| \le \liminf_{n\to\infty}\|A_{n}^{-1}\|\,.
  \end{equation}

  In particular, if $p = 2$ and $A$ and $\fol{A}{n}$ are self-adjoint operators on $\ell^2(\II)$ with $A_n x \to A x$ for all $x\in \cfin(\II)$, then the above conclusion remains valid.
\end{lemma}
\begin{proof}
We follow the proof in~\cite{{Lindner.2006}} but also involve Lemma~\ref{lem:lowerNorm}.

Let $(A_{n_k})_{k \in \NN}$ be a stable subsequence such that 
\begin{equation*}
  \lim_{k \to \infty} \| A_{n_k}^{-1} \| = \liminf_{n \to \infty} \|A_n^{-1}\|\, . 
\end{equation*}
This gives that, for all $x \in \cfin(\II)$ and all $k \geq k_0$ such that $A_{k}$ is invertible,
\begin{equation*}
  \|x\| 
  = \|A_{n_k}^{-1}A_{n_k}^{} x\| 
  \le \|A_{n_k}^{-1}\| \|A_{n_k}^{}x\| 
  \to \liminf_{n \to \infty} \|A_n^{-1}\| \|Ax\| \quad\text{for } k \to \infty.
\end{equation*}
In particular, we have
\begin{equation*}
\Big(\liminf_{n \to \infty} \|A_n^{-1}\| \Big)^{-1} \leq \nu_D(A), \quad\text{for all } D \in \NN.
\end{equation*}
From Lemma~\ref{lem:lowerNorm}, we deduce that,
\begin{equation*}
\Big(\liminf_{n \to \infty} \|A_n^{-1}\| \Big)^{-1} \leq \nu(A) + \varepsilon, \quad\text{for all } \varepsilon > 0,
\end{equation*}
and, consequently,
\begin{equation*}
\Big(\liminf_{n \to \infty} \|A_n^{-1}\| \Big)^{-1} \leq \nu(A).
\end{equation*}
Similarly, we obtain~$\big(\liminf_{n \to \infty} \|A_n^{-1}\| \big)^{-1} \leq \nu(A')$.
Now, Lemma~\ref{lem:lowerNormInvertible} allows us to conclude~\eqref{eq: lim inf norm inverse estimate}
\end{proof}

\section{Basic Properties of Sturmian Words}\label{sec:sturm}
In order to study discrete aperiodic Schrödinger operators $H \colon \ell^p(\ZZ) \rightarrow \ell^p(\ZZ)$ with
\begin{equation*}
	(H x)_n = x_{n+1} +   x_{n-1} 	 + v(n) x_n \, ,\qquad n\in\ZZ \,,
\end{equation*}
we start by introducing the class of Sturmian potentials $v$. For this we leave the operator theoretic setting and only consider the potential as a two-sided infinite word over the alphabet $\Sigma=\{0,1\}$. The combinatorial properties, proved in this section, will be the groundwork for our main results.

For $\alpha \in [0,1]$ and $\theta \in [0,1)$, we consider, as in~\eqref{eq:defintion_words}, the two-sided infinite words $v_{\alpha,\theta}$ and $\widetilde{v}_{\alpha,\theta}$ defined by
\begin{equation}\label{eq:vdef}
v_{\alpha,\theta} (k) \coloneqq \chi_{[1-\alpha,1)}(k\alpha+\theta \bmod 1)
\end{equation}
and 
\begin{equation}\label{eq:vdefalt}
  \widetilde{v}_{\alpha,\theta} (k) \coloneqq \chi_{(1-\alpha,1]}(k\alpha+\theta \bmod 1)
\end{equation}
for $k\in \ZZ$. 

For irrational~$\alpha$, the words $v_{\alpha,\theta}$ and $\widetilde{v}_{\alpha, \theta}$ are called \emph{Sturmian words} and $\theta$ is called the associated offset.
Observe that these two words satisfy a symmetry property, namely, for all $k\in\ZZ$, we have
\begin{align}\label{eq:mirroredLimWord}
v_{\alpha, \theta}(k)
  &= \chi_{[1-\alpha,1)}(k\alpha+\theta\bmod 1)\notag \\
  &=\chi_{(1-\alpha,1]}(-k\alpha- \alpha-\theta \bmod 1) \\
  &=\widetilde{v}_{\alpha,(-\alpha-\theta\bmod 1)}(-k)\,.\notag
\end{align}
If $\theta = 0$, we use the shorthand notation $v_\alpha \coloneqq v_{\alpha,0}$ and $\widetilde{v}_{\alpha} \coloneqq \widetilde{v}_{\alpha,0}$. 

For irrational~$\alpha$, we consider its \emph{continued fraction expansion}, cf.~\cite{Khinchin.1997},
\begin{equation*}
\alpha = \frac{1}{a_1 + \frac{1}{a_2 + \frac{1}{a_3 + \dots}}}  \eqqcolon [a_1,a_2,a_3,\dots]
\quad\text{with}\quad
a_1,a_2,\ldots\in\NN.
\end{equation*}	
If we terminate this infinite fraction expansion for a given $m \in \NN$, \ie we remove the term $+\frac 1{a_{m+1}+\dots}$, 
we get a rational number $\alpha_m\coloneqq \frac{p_m}{q_m}$ which is called the \emph{$m$-th rational approximant} of $\alpha$. It satisfies the following recursion
	\begin{align*}
	&p_{-1}\coloneqq  1 \,, ~p_0\coloneqq  0 \,, ~ p_m \coloneqq a_m p_{m-1} + p_{m-2}\, , \\
	&q_{-1}\coloneqq  0 \,, ~q_0\coloneqq  1 \,, ~ q_m \coloneqq a_m q_{m-1}+q_{m-2}  \,.
	\end{align*}
One can show, see~\cite[Theorem~5]{Lang.1995} or~\cite[Theorem 9, Theorem 13]{Khinchin.1997}, that the approximation of $\alpha$ with the rational numbers $\frac{p_m}{q_m}$
	satisfies
	\begin{equation}\label{eq:approximation_rate1}
        \frac{1}{q_m (q_{m + 1} + q_m)} < \Big| \alpha - \frac{p_m}{q_m} \Big| <  \frac{1}{q_m q_{m+1}} \leq \frac{1}{m(m+1)}
	\end{equation}
	and also
        \begin{equation}\label{eq:approximation_rate2}
		\frac{p_0}{q_0} < \frac{p_2}{q_2} < \cdots < \alpha < \cdots < \frac{p_3}{q_3} < \frac{p_1}{q_1} \,.
        \end{equation}

The rational approximant $\alpha_m\coloneqq \frac{p_m}{q_m}$ is a best approximation, in the sense that              
\begin{equation*}
			\| q_m \alpha \| = | q_m \alpha - p_m | 
			\quad \text{and} \quad
			\| q \alpha \| > 	\| q_m \alpha \|
			~
			\text{ for }
			1 \leq q< q_m \,,
\end{equation*}
where $\| x \| \coloneqq \inf \{|x - k| : k \in \ZZ\}$ denotes the distance of $x$ to $\ZZ$.

\begin{example}\label{ex:goldenRatio}
  For $\alpha = (\sqrt{5} - 1)/2$, the corresponding sequence of $m$-th rational approximants yields
  \begin{equation*}
    0 
    < \frac{1}{2} < \frac{3}{5} < \frac{8}{13} < \frac{21}{34} 
    < \cdots 
    < \alpha 
    < \cdots 
    < \frac{34}{55}< \frac{13}{21} < \frac{5}{8}< \frac{2}{3} < 1 \,.
  \end{equation*}
\end{example}

Let us summarise these notations and assumptions in the following hypothesis:
\begin{hypo}\label{hypo:words}
	Let $\alpha \in [0,1]$ irrational, $\theta \in [0,1)$, and let $\alpha_m = \frac{p_m}{q_m}$ denote the $m$-th rational approximant to $\alpha$. In addition, let $v_{\alpha,\theta}$ and $\widetilde{v}_{\alpha, \theta}$ be Sturmian words as in \eqref{eq:vdef} and \eqref{eq:vdefalt}.
	Further, the periodic words $v_{\alpha_m,\theta} $ and $\widetilde{v}_{\alpha_m,\theta} $ are defined by
	\begin{align}
	v_{\alpha_m,\theta} (k) 
  &\coloneqq \chi_{[1-\alpha,1)}(k\alpha_m+\theta \bmod 1)	\quad
	\textnormal{and} \\
	\widetilde{v}_{\alpha_m,\theta} (k) 
  &\coloneqq \chi_{(1-\alpha,1]}(k\alpha_m+\theta \bmod 1)\,.
	\end{align}
\end{hypo}        

In Section~\ref{sec:rationalApproximations}, we use the properties of the rational approximants $\alpha_m$ to relate the words $v_{\alpha_m,\theta}$ and $v_{\alpha,\theta}$ for particular values of $\theta$.
In Section~\ref{sec:limitWords}, we consider subshifts generated by the word $v_\alpha$ and their finite subwords.

\subsection{Rational Approximations}\label{sec:rationalApproximations}
Given an irrational $\alpha \in [0,1]$ with corresponding continued fraction expansion $[a_1, a_2, a_3, \dots]$, we define words $s_m$ of finite length over the alphabet $\Sigma=\{0,1\}$ by the recursion:
\begin{equation}\label{eq:recursiveWord}
	s_{-1} \coloneqq 1,\; s_0\coloneqq 0,\; s_1 \coloneqq s_0^{a_1-1}s_{-1},\; s_m\coloneqq s_{m-1}^{a_m}s_{m-2}.
\end{equation}
Note that, for all $m\geq 0$, the length of the word $s_m$ is equal to $q_m$, where $\frac{p_m}{q_m}$ is the $m$-th rational approximant to $\alpha$. 

\begin{prop}[{\cite[Proposition 4.5]{Berstel.1996}}]\label{prop:Palindrom}
Let $\alpha \in [0,1]$ be irrational and consider the words $s_m$ defined by~\eqref{eq:recursiveWord}. 
  Then there exist palindromes $\pi_m$, $m\geq 2$, such that
	\begin{equation*}
	s_{m} = \begin{cases} \pi_m 1 0 \, , & m  \text{ even,} \\  
        \pi_m 0 1 \, , & m \text{ odd.} \end{cases}
	\end{equation*} 
\end{prop}
For $m\geq 2$, the word $s_{m-1}$ is a prefix of $s_m$. Therefore the following limit from the left 
exists in an obvious sense:
	\begin{align*}
	v_+ &= \lim_{m \to \infty} s_m.
	\end{align*} 	
	
The finite words $s_m$ are in particular important, since, for the offset $\theta = 0$, the Sturmian word $v_{\alpha,0}$, restricted to $\{1,2,3,\dots\}$, coincides with the one-sided infinite word $v_+$, see \cite{Bellissard.1989}. Due to this equality and the recursive definition of $s_m$, Sturmian words are locally periodic. This, together with a symmetry property, is stated in the following lemma by Bellissard et al. in \cite{Bellissard.1989}.
Note that the original result is only formulated for $v_{\alpha}$ but the proof for $\widetilde{v}_{\alpha}$ is analogous. 

\begin{lemma}\label{lem:Bellissard01}
Let $\alpha \in [0,1]$ be irrational, $m\geq 2$. Then, for $v\in \{v_{\alpha},\widetilde{v}_{\alpha}\}$, the following holds:
\begin{enumerate}[label=(\alph*)]
                \item\label{it:bellissard1b} 
                For all $1 \leq k \leq q_{m + 1} - 2$, we have $v(q_m + k) = v(k)$.
        \item\label{it:extensionneu} For all $k\geq 2$,  we have $v (-k) = v(k-1)$.  
            \end{enumerate}
\end{lemma}	

Due to the symmetry property \ref{it:extensionneu} in Lemma \ref{lem:Bellissard01}, and since one can decompose the words $s_m$ using a palindrome, we obtain that the Sturmian word $v_{\alpha,0}$, restricted to $\{\dots,-2,-1,0\}$, is equal to the one-sided infinite word $v_-$ given by the following limit from the right
	\begin{align*}
	v_- &\coloneqq \lim_{m \to \infty} s_{2m}.
	\end{align*} 	
Note that this limit exists since the  word $s_{2m-2}$ is a suffix of $s_{2m}$. 

The following result will be useful to prove the next lemmata in this section. 
Note that here we do not assume $\alpha$ to be irrational unless explicitly stated.
	\begin{lemma}\label{lem:modcalculation} 
	 Let $\alpha \in [0,1]$, $k\in \ZZ$, and $|\theta| < \frac{1}{2}$. We have 
          \begin{equation*}
            v_{\alpha,0}(k) = v_{\alpha,\pm \theta \bmod 1}(k)\quad\text{and}\quad \widetilde{v}_{\alpha,0}(k) = \widetilde{v}_{\alpha,\pm \theta \bmod 1}(k)
          \end{equation*}
          if one of the following conditions is fulfilled:
\begin{enumerate}[label=(\alph*)]
                      \item\label{it:gtgt} $\norm{k \alpha} > |\theta|$ and $\norm{(k+1) \alpha} > |\theta|$.
                      \item\label{it:geqgt} $\norm{k \alpha} \geq |\theta|$, $\norm{(k+1) \alpha} >  |\theta|$,
	 		and $k \alpha \pm \theta \notin \ZZ$.
                      \item\label{it:gtgeq} $\norm{k \alpha} > |\theta|$, $\norm{(k+1) \alpha} \geq  |\theta|$,
	 		and $(k+1) \alpha \pm \theta \notin \ZZ$.
	 	\end{enumerate}	
	\end{lemma}	
	
	\begin{proof}
	We only show the result for $v_{\alpha,\theta}$ because, for $\widetilde{v}_{\alpha,\theta}$, the proof works completely analogously.
  For part~\ref{it:gtgt}, let $|\theta| < \frac{1}{2}$ and $k\in \ZZ$ be such that $\norm{k \alpha} > |\theta|$ and $\norm{(k+1) \alpha} > |\theta|$. 
          Then we have the following equivalences:
		\begin{align*}
          k \alpha &\bmod 1  \in [1-\alpha,1) \\
					&~\iff~      	\exists m \in \ZZ \text{ with }
          1-\alpha  \leq k\alpha - m < 1   \\
            	&~\iff~      	\exists m \in \ZZ \text{ with }
          k \alpha   < m+1 \leq (k+1) \alpha   \\
          	&~\iff~ \exists m \in \ZZ \text{ with }
          k \alpha \pm \theta < m+1 \leq (k+1) \alpha \pm \theta \\
			&~\iff~  k \alpha \pm \theta \bmod 1  \in [1-\alpha,1)   
          \,.
		\end{align*}
Part~\ref{it:geqgt} and~\ref{it:gtgeq} follow from part~\ref{it:gtgt} by using the fact that, from $k\alpha \pm \theta \notin \ZZ$ or $(k+1)\alpha \pm \theta \notin \ZZ$, 
we conclude $\|k\alpha\| \neq \pm \theta$ or $\|(k+1)\alpha\| \neq \pm \theta$, respectively.
	\end{proof}

        The following lemma shows that the \emph{rational approximation} of irrational numbers $\alpha \in [0,1]$ via the sequence $(\alpha_n)$, gives rise to a \emph{periodic approximation} of Sturmian potentials $v_{\alpha, 0}$ via periodic words $v_{\alpha_n, 0}$.
		
	\begin{lemma}\label{lem:potential_aperiodic1}   
	Let $\alpha\in [0,1]$ be irrational, $n\geq 2$, and $\alpha_n = \frac{p_n}{q_n}$ the $n$-th rational approximant to $\alpha$. Then
                \begin{equation*}
	 	 v_{\alpha_n}(k) = v_{\alpha}(k) \quad \text{ for all } 
	 	 \begin{cases}-q_n + 1 \leq k \leq q_{n+1}-2  \, , & n  \text{ even,} \\  -q_{n+1}-1 \leq k \leq q_n - 2 \, , & n \text{ odd,} \end{cases}   	
                \end{equation*} 
and 					
          \begin{equation*}
	 	 \widetilde{v}_{\alpha_n}(k) = \widetilde{v}_{\alpha}(k) \quad \text{ for all } 
	 	 \begin{cases} -q_{n+1} -1 \leq k \leq q_n - 2 \, , & n \text{ even,} \\
		-q_n + 1 \leq k \leq q_{n+1}-2  \, , & n  \text{ odd.} \\
		\end{cases}
          \end{equation*}			
	\end{lemma}
	
	\begin{proof}	
	We first show that $v_{\alpha_n}(k) = v_{\alpha}(k)$ for $1\le k \le q_n-2$. 
        To this end, let us fix $1\le k \le q_n-2$ and set $\theta_k \coloneqq k(\alpha-\alpha_n)$. 
          Using the estimates~\eqref{eq:approximation_rate1}, yields 
          \begin{equation*}
            |\theta_k| < \frac{1}{2} \quad\text{and}\quad \min\big\{\|k\alpha_n\|, \|(k+1)\alpha_n\|\big\} \geq \frac{1}{q_n} > |\theta_k|\,. 
          \end{equation*}
          We now use Lemma~\ref{lem:modcalculation}\ref{it:gtgt} to derive the identity
          \begin{equation}\label{eq:smallcoincidence}
            v_{\alpha_n,0}(k)=v_{\alpha_n,\theta_k \bmod 1}(k)= v_{\alpha,0}(k) \quad\text{for all } 1\le k \le q_n-2\,.
          \end{equation}
	
          By using the recurrence properties of $v_{\alpha}$, we now extend the validity of~\eqref{eq:smallcoincidence} to a larger set. 
Recall that, by Proposition~\ref{prop:Palindrom}, there is a palindrome~$\pi_n$ such that $s_n = \pi_n 10$ if $n$ is even and $s_n = \pi_n 01$ if $n$ is odd. 
Since $v_{\alpha_n}(0) = v_{\alpha}(0)=0$ and $v_{\alpha_n}(-1) = v_{\alpha}(-1)=1$ and since both words coincide for $1\le k \le q_n-2$, we get that the period of $v_{\alpha_n}$ is equal to $\pi_n 10$. 
Now, using Lemma~\ref{lem:Bellissard01} for $v_\alpha$, we get for $n$ even the following diagramme:
\begin{center}
	{\def\arraystretch{1.5}\tabcolsep=1.25ex
		\begin{tabular}{c@{\hskip 3ex}c@{\hskip 3ex}c@{\hskip 3ex}c@{\hskip 3ex}c@{\hskip 3ex}c@{\hskip 3ex}c@{\hskip 3ex}c}
			\toprule
			\rowcolor{gray!45}
                        $k$ &  $-1-q_n \ldots$  &  $-1$ & $0$ & $1 \ldots q_n$  &$\ldots$   & $\dots a_{n + 1} q_n$ & $\ldots q_{n + 1}$ 
			\\
			\bottomrule
                        $v_\alpha(k)$  & $s_n^{\mathrm{R}}$ &  $ \hphantom{-}1$ & $0$  & $s_n$  & $\ldots$ & $ s_n$ & $s_{n-1}$ 
			\\ \rowcolor{gray!25}
                         $v_\alpha(k)$ & $0~\mathbf{1}~~~\pi_n$&  $ \hphantom{-}1$ & $0$  & $\pi_n~~1~0$  &  $\ldots$ & $\pi_n~~1~0$ & $\pi_{n-1}~~\mathbf{0}~1$ 
			\\
			[0.1cm]
			\bottomrule
		\end{tabular}
	}
\end{center}
Here, $s_n^{\mathrm{R}}$ denotes the reversed word of $s_n$.
The bold numbers mark the first positions where $v_\alpha$ does not coincide with the periodic word $v_{\alpha_n}$. 
For $n$ odd, we have $s_n = \pi_n 01$ and the period of $v_{\alpha_n}$ is now equal to $s_n^{\mathrm{R}}$. 
Hence, $v_{\alpha_n}(k) = v_{\alpha}(k)$ for $-q_{n+1}-1 \leq k \leq q_n - 2$. 
This proves the first part of the lemma.

The proof for $\widetilde{v}_{\alpha_n}$ and $\widetilde{v}_{\alpha}$ works similarly. 
Here, we have to use that $\widetilde{v}_{\alpha}(k)$ and $v_{\alpha}(k)$ only differ for $k\in \{-1,0\}$ and that the period of $\widetilde{v}_{\alpha_n}$ now is equal to $\pi_n 01$.
\end{proof}

\begin{remark}\label{rem:shift_invariant}
Note that, by the recursion $s_n=s_{n-1}^{a_n}s_{n-2}$ and the above reasoning, we have $v_\alpha(k) = v_{\alpha}(k-q_n)$ for $n$ even and all $1\le k \le q_{n+1}$.
\end{remark}

Combining Lemma~\ref{lem:modcalculation} and Lemma~\ref{lem:potential_aperiodic1}, we can show the following:	
	
\begin{lemma}\label{lem:potential_aperiodic_theta} 
Assume Hypothesis~\ref{hypo:words} with $m\geq 2$.  

\begin{enumerate}[label=(\alph*)] 
 \item\label{it:kgt0} If  $\theta = k_0\alpha  \bmod 1$ for some $k_0\in\ZZ_-$ choose $m$ sufficiently large such that  $ |k_0| \le q_{m-1}$.
Then for $-q_m+1 \le k \le q_m-2$, the following holds:

 \begin{equation*}
v_{\alpha_m,\theta}(k) = \begin{cases} \widetilde{v}_{\alpha,\theta}(k)   &  \text{if $m$ is even,} \\ 
  v_{\alpha_m,\theta}(k) = {v}_{\alpha,\theta}(k)& \text{ if $m$ is odd.} \end{cases}  
 \end{equation*}

  \item\label{it:notorbit} If $\theta\not\in\{k_0\alpha \bmod 1 : k_0\in\ZZ_-\}$ then for each $D\in\NN$  there is an $m_0$ such that for all $m\ge m_0$ and for all $k=0,\dots, D$
 \begin{equation*}v_{\alpha_m,\theta}(k) = v_{\alpha,\theta}(k). \end{equation*}
\end{enumerate}
\end{lemma}	
\begin{proof}
  Let $k_0 \in \ZZ$ with $q_{m - 1} \geq |k_0|$, and set $\theta_m \coloneqq k_0 \alpha_m \bmod 1$.
  Then, for all $k \in \ZZ$,
  \begin{equation*}
    v_{\alpha, \theta} (k) =v_{\alpha, 0}(k+k_0) \quad\text{and}\quad v_{\alpha_m,\theta_n } (k) =v_{\alpha_m, 0}(k+k_0)\,,
  \end{equation*}
and, by Lemma~\ref{lem:potential_aperiodic1}, we conclude
  \begin{align}\label{eq:alpha0coincidence}
    v_{\alpha, 0} (k + k_0) =  v_{\alpha_m, 0} (k + k_0) \quad\text{for} -q_n+1  \le k+k_0\le q_n-2 \,. 
  \end{align}
  As $q_n \geq n \geq 2$ by assumption, identity~\eqref{eq:alpha0coincidence} implies that the $q_n$-periodic word $v_{\alpha_m,0}$ coincides with $v_{\alpha, 0}$ for at least one full period.
  This fact and Lemma~\ref{lem:Bellissard01} allow us to extend the region of indices specified in~\eqref{eq:alpha0coincidence} to the left and right by one further period  such that the assumption $q_{n - 1} \geq |k_0|$ implies
  \begin{equation*}
    v_{\alpha, \theta} (k) =  v_{\alpha_m, \theta_n} (k) \quad\text{for} -q_n+1  \le k\le q_n-2 \,. 
  \end{equation*}
  Similar to the proof of Lemma~\ref{lem:potential_aperiodic1}, we use a perturbation argument to pass from $v_{\alpha_m,\theta_n}$ to $v_{\alpha_m,\theta}$. We consider two cases.

\emph{Case 1: assume that $(k+k_0)\alpha_m \not\in \ZZ$ and $(k+k_0+1)\alpha_m \not\in \ZZ$.}
We define $x\coloneqq k_0(\alpha-\alpha_m)$. 
Then the approximation rate estimates~\eqref{eq:approximation_rate1} give that
  \begin{equation*}
    |x| < \frac{1}{q_n}\le \min \big\{\|(k_0+k)\alpha_m\|,\|(k_0+k+1)\alpha_m \|\big\}.
  \end{equation*}
  Now, Lemma~\ref{lem:modcalculation}\ref{it:gtgt} yields
  \begin{equation*}
    v_{\alpha_m,\theta_n}(k) = v_{\alpha_m, 0}(k+k_0)= v_{\alpha_m,x \bmod 1}(k+k_0) = v_{\alpha_m, \theta}(k) .
  \end{equation*}
  This shows
  \begin{equation*}
    v_{\alpha, \theta} (k) =  v_{\alpha_m, \theta} (k) \quad\text{for} -q_n+1  \le k\le q_n-2 \,. 
  \end{equation*}
  Note that, for this range of $k$, we also have $v_{\alpha, \theta} (k) =\widetilde{v}_{\alpha, \theta} (k) $. 

\emph{Case 2: assume that $(k+k_0)\alpha_m \in \ZZ$ or $(k+k_0+1)\alpha_m \in \ZZ$.}
 We first check the case $k+k_0 = j q_m$ for $j \in\ZZ$. 
 Then $v_{\alpha_m, \theta_m} = 0$ and  $\widetilde{v}_{\alpha_m, \theta_m} = 1$. 
 We get that
 \begin{align*}
  v_{\alpha_m,\theta}(k) &= \chi_{[1-\alpha_m, 1)}((jq_m-k_0)\alpha_m+k_0 \alpha \bmod 1) 
\\&= \chi_{[1-\alpha_m, 1)}(k_0 (\alpha-\alpha_m) \bmod 1) 
\\& = \begin{cases}
            \chi_{[-\alpha_m, 0)}(k_0 (\alpha-\alpha_m) )&=    1 =\widetilde{v}_{\alpha_m, \theta_m}(k) \quad m \text{ even,}  \\ 
                            \chi_{[-\alpha_m, 0) \cup [1-\alpha_m, 1)}(k_0 (\alpha-\alpha_m) ) &=    0 =  {{v}_{\alpha_m, \theta_m}}(k) \quad m \text{ odd,}
                                           \end{cases} 
 \end{align*}
holds. For $k+k_0+1 = j q_m$ for $j \in\ZZ$, we compute
 \begin{align*}
  v_{\alpha_m,\theta}(k) &= \chi_{[1-\alpha_m, 1)}((jq_m-k_0-1)\alpha_m+k_0 \alpha \bmod 1) 
\\&= \chi_{[1-\alpha_m, 1)}(k_0 (\alpha-\alpha_m) -\alpha_m \bmod 1) 
\\& = \begin{cases}
   \chi_{[-\alpha_m, 0)}(k_0 (\alpha-\alpha_m)-\alpha_m ) &= 0 =  \widetilde{v}_{\alpha_m, \theta_m}(k) \quad m \text{ even,} \\ 
  \chi_{[-\alpha_m, 0) \cup [1-\alpha_m, 1)}(k_0 (\alpha-\alpha_m)-\alpha_m ) &= 1 =  {v}_{\alpha_m, \theta_m}(k) \quad m \text{ odd.} \end{cases}                             
 \end{align*}
 This proves the case~\ref{it:kgt0}.

Finally, we prove \ref{it:notorbit}, i.e. the case that  $\theta\not\in\{k_0\alpha \bmod 1 : k_0\in\ZZ_-\}$. 
Let $D\in\NN$ and set $\varepsilon \coloneqq 	\min_{0\le k\le D} \| k \alpha +\theta\| >0$. 
Choose $n_0$ sufficiently large such that $\|(D+1)(\alpha-\alpha_{n_0})\| < \varepsilon$.
Then, for $n\ge n_0$ and $0\le k\le D$,
\begin{align*}
 v_{\alpha, \theta}(k)&=\chi_{[1-\alpha, 1)}(k \alpha + \theta \bmod 1)   \\
 &=\chi_{[1-\alpha, 1)}(k \alpha+ \theta  - (k+1)(\alpha-\alpha_n) \bmod 1)   \\
  &=\chi_{[1-\alpha, 1)}(k \alpha_n+ \theta  -\alpha + \alpha_m \bmod 1)   \\
 &= v_{\alpha_{n},\theta}(k). \qedhere
\end{align*}
\end{proof}

\subsection{Subwords and Limit Words}\label{sec:limitWords}   
In this subsection, we analyse subshifts of the potential similar to~\cite{Lindner.2018} and their finite subwords.
The observations we make will translate to a characterisation of limit operators of the Schr\"odinger operator with Sturmian potential and their spectra, see Section~\ref{sec:limOps}.

	Let $v\coloneqq (v_{\alpha,\theta} (n))_{n\in \ZZ}$ be the Sturmian potential. Then we can write
        \begin{equation*}
	H_{\lambda,\alpha,\theta} = S^{-1} + S  + \lambda \, M_v
	\, , \quad x \in \ell^p(\ZZ) \, ,
        \end{equation*}
    where $S$ denotes the right shift operator, i.e.\ $(Sx)(n) = x(n - 1)$, and $M_b$ denotes the multiplication operator by $b\in \ell^\infty(\ZZ)$.      
Let $h = (h_1, h_2, \dots)$ be a sequence in $\ZZ$ with $h_k \to \pm \infty$, so that the limit operator $(H_{\lambda,\alpha,\theta})_h$ of $H_{\lambda,\alpha,\theta}$ exists. Then
\begin{equation*}
(H_{\lambda,\alpha,\theta})_h = S^{-1} + S + \lambda \, M_{v_h}
\end{equation*}
with a new potential
\begin{equation*}
  v_h \coloneqq  \lim_{k\to\infty} S^{-h_k} v
\end{equation*}
where the limit is taken w.r.t.\ pointwise convergence.
The set 
\begin{equation}
\cF_v \coloneqq  \Big\{\lim_{k\to\infty} S^{-h_{k}} v : h_k \to \pm \infty\Big\}
\end{equation}
of all such potentials $v_h$ is the so-called \emph{subshift generated by the sequence $v$}.
The same can be done for periodic words $v=v_{\alpha_m, \theta}$. In that case, the subshift generated by $v$ is finite. 

For Sturmian potentials, the subshift is explicitly known, see, e.g.,~\cite[Theorem~2.14]{Damanik.2007} and the appendix of~\cite{DamanikLenz.2004}.

\begin{proposition}\label{prop:subshift}
Assume Hypothesis~\ref{hypo:words}. Then we have the following:
\begin{enumerate}[label=(\roman*)]
 \item 
    The subshift of $v_{\alpha, \theta}$ is given by
\begin{equation}\label{eq:subshifts}
  \cF_{v_{\alpha,\theta}}  = \big\{v_{\alpha,\psi}, \widetilde{v}_{\alpha,\psi} : \psi\in [0,1)\big\}\,.
\end{equation}
    In particular, $\cF_{v_{\alpha,\theta}} = \cF_{v_{\alpha,0}} = \cF_{\widetilde{v}_{\alpha,\theta}}$.
\item
    The subshift of $v_{\alpha_m, \theta}$ is given by
\begin{equation}\label{eq:subshiftsperiodic}
  \cF_{v_{\alpha_m,\theta}}  = \big\{v_{\alpha_m,\psi}, \widetilde{v}_{\alpha_m,\psi} : \psi\in [0,1)\big\}\,.
\end{equation}
    In particular, $\cF_{v_{\alpha_m,\theta}} = \cF_{v_{\alpha_m,0}} = \cF_{\widetilde{v}_{\alpha_m,\theta}}$.
\end{enumerate}

\end{proposition}

As equation~\eqref{eq:subshifts} shows, the subshift $\cF_{v_{\alpha, \theta}}$ does not depend on $\theta$ and also does not change when switching to the potential $\widetilde{v}_{\alpha,\theta}$.
Therefore, the shorthand notation $\cF_\alpha = \cF_{v_{\alpha,\theta}}$ is well defined.
The same is true for $\cF_{\alpha_m} = \cF_{v_{\alpha_m,\theta}}$.

Next, we examine the finite subwords of words in the subshifts $\cF_{\alpha}$ and $\cF_{\alpha_m}$.
For a word $v$ over $\ZZ$, we define $v_+$ to be the restriction of $v$ onto $\NN$ and $v_-$ to be the restriction onto $\ZZ\setminus \NN$.
We write $w\prec v$ when $w$ is a subword of $v$,  cf.~\cite[Section 3.1]{Lindner.2018}.

 \begin{lemma}\label{lem:subwordInRightHalf}
 Assume Hypothesis~\ref{hypo:words}.
 Let $\cF$ be either $\cF_{\alpha}$ or $\cF_{\alpha_m}$. 
 Then for $v\in \cF$ and $w \prec v$ we have that
 
 \begin{enumerate}[label=(\roman*)]
  \item \label{it:subwords} $w \prec u$ for all $u\in \cF$ and $w$ occurs infinitely many times in $u$,
  \item \label{it:subwordright} $w \prec v_+$ and $w$ occurs infinitely many times in $v_+$,
  \item \label{it:subwordleft} $w \prec v_-$ and $w$ occurs infinitely many times in $v_-$.
 \end{enumerate}
\end{lemma}
\begin{proof}
Part~\ref{it:subwords} follows from {\cite[Proposition 2.1.18]{Lothaire.2002}}. 
For part~\ref{it:subwordright}, take $w\prec {v}_{\alpha, \theta}$ and pick $l,r\in\ZZ$ such that $\supp(w) \subset [l,r]$. 
Choose $k_n\subset\NN$ such that $k_n\alpha +\psi \bmod 1$ converges to $ l\alpha +\theta \bmod 1 $ from above. 
Since $\chi_{[1-\alpha,1)}$ is right-continuous, there is an $n_0\in\NN$ such that for all $n\ge n_0, j\in[0,r-l]$ we have that
\begin{align*} \begin{aligned}
  v_{\alpha, \psi}(k_n+j)&=\chi_{[1-\alpha,1)}((k_n+j)\alpha +\psi  \bmod 1) \\
  &= \chi_{[1-\alpha,1)}((l+j)\alpha +\theta  \bmod 1) = v_{\alpha, \theta}(l+j).
\end{aligned}\end{align*}
Hence, $w\prec (v_{\alpha, \psi})_+$.
The proof for \ref{it:subwordleft} is analogous.
\end{proof}

For an infinite word $v$, let $\mathcal{W}_D(v)$ denote the set of all subwords of $v$  that have length $D\in\NN$.
\begin{lemma}\label{lem:sameSubwordsApprox}
    Assume Hypothesis~\ref{hypo:words}.
    Let $v  \in\cF_\alpha$ and $v_m \in \cF_{\alpha_m}$.
    Then, for each $D \in \NN$, there exists $N \in \NN$, such that, for all $m \geq N$, we have $\mathcal{W}_D(v) = \mathcal{W}_{D}(v_m)$.
\end{lemma}

\begin{proof}
We only prove the case $v= v_{\alpha,0}$ since the other cases follow by Lemma~\ref{lem:subwordInRightHalf}.
    Let $D \in \NN$ and choose $m \in \NN$ such that $q_m \geq 2D$ and $\mathcal{W}_D(v) \subset \mathcal{W}_D(v_m)$.
    Note that this is possible since $v$ has precisely $D + 1$ subwords of length $D$, see~\cite{Berstel.1996}.
    Now let $w \in \mathcal{W}_D(v_m)$.
    Since $v_m$ is periodic, we may assume that 
    \begin{align*}
        w \prec (v_m(-q_m + 1), \dots, v_m(q_m - 2))\;.
    \end{align*}
    Now the claim follows from Lemma~\ref{lem:potential_aperiodic1}.
\end{proof}

\begin{lemma}\label{lem:thetaConvergence}
 Assume Hypothesis~\ref{hypo:words}. For $n\in\NN$, let  $v_n\in \cF_\alpha $. Let further $\theta_n$ denote the offset of $v_n$.
 \begin{enumerate}[label=(\roman*)]
  \item \label{it:orbitTF1} Then there is a subsequence $(v_{n_j})$ that converges pointwise, i.e.\ $v_{n_k}(i)$ converges for every $i\in \ZZ$.
  \item \label{it:orbitTF2} If $\{\theta_n : n\in\NN\}$ is dense in $[0,1]$, then, for each $v\in\cF_\alpha$, there is a subsequence $(v_{n_j})$ that converges pointwise to $v$.
 \end{enumerate}
\end{lemma}
\begin{proof}
To prove part~\ref{it:orbitTF1}, set $v_n = v_{\alpha, \theta_n}$ or $v_n= \widetilde{v}_{\alpha,\theta_n}$ with $\theta_n\in [0,1)$ for all $n\in \NN$.
Now, pass to a subsequence $v_{n_j}$ such that $\theta_{n_j}$ converges to some $\psi\in[0,1]$ either from the left or from the right.

If there is a constant subsequence in $\theta_{n_j}$ then the assertion is immediate.
For the other case, fix $k\in\ZZ$ and take $j_0$ sufficiently large such that 
\begin{equation}
k\alpha-\theta_{n_j}  \bmod 1 \not \in \{0,1-\alpha\} \quad \textnormal{for}  \quad j> j_0 \,.
\end{equation}
If the convergence $\theta_{n_j}\rightarrow \psi$ is from the right, then, due to the right-continuity of $\chi_{[1-\alpha,1)}(k\alpha- \cdot \bmod 1)$, we find that 
\begin{equation}
v_{n_j}(k) = \chi_{[1-\alpha,1)}(k\alpha- \theta_{n_j} \bmod 1) \rightarrow  \chi_{[1-\alpha,1)}(k\alpha- \psi \bmod 1) = v_{\alpha, \psi}(k) \,. 
\end{equation}

If the convergence is from the left and $\psi < 1$, then we find that $v_{n_j}(k)\rightarrow \widetilde{v}_{\alpha, \psi}(k)$, using the left-continuity of $\chi_{(1-\alpha,1]}(k\alpha- \cdot \bmod 1)$.
If $\psi = 1$, then $v_{n_j}(k)\rightarrow \widetilde{v}_{\alpha, 0}(k)$ with the same ideas.

For part~\ref{it:orbitTF2}, let $\psi\in[0,1)$ be the offset of $v$.
If $v=v_{\alpha, \psi}$, we choose a subsequence $(\theta_{n_j})$ that converges to $\psi$ from the right and argue as in the proof of part~\ref{it:orbitTF1} that $v_{n_j} \rightarrow v$ pointwise.
For the case that $v=\widetilde{v}_{\alpha, \psi}$, we do the same with convergence from the left. 
If $\psi=0$, then we take a subsequence $(\theta_{n_j})$ that converges to $1$.
\end{proof}

\section{Spectrum of \texorpdfstring{$H_{\lambda, \alpha,\theta}$}{H}}\label{sec:specHltheta}

In this section, we further develop the framework to be used in order to characterise the applicability of the FSM.
Recall that the generic FSM described in Section~\ref{sec:fsm_intro} is an approximation method resembling the classical strategy of using finite rank operators to approximate an infinite rank operator~$H_{\lambda, \alpha, \theta}$.
In order to tackle the invertibility problem of limit operators that lies at the core of the FSM, see Proposition~\ref{prop:FSMCondition}, we will work with a different approximation method which, nevertheless, also comes natural in our setting:
the approximation of aperiodic Schrödinger operators via periodic operators of infinite rank.
This method resembles the operator theoretic consequence of the periodic approximation $v_{\alpha_m,\theta}$ of Sturmian words $v_{\alpha,\theta}$ from Section~\ref{sec:rationalApproximations}.

There are three types of results in this section:
We start by characterising the limit operators of aperiodic Schrödinger operators $H_{\lambda, \alpha, \theta}$ and examine the spectral quantities of them and their one-sided compressions.  
We follow this up with Section~\ref{sec:traceCondition}, where we give a condition for an energy~$E$ to lie in the resolvent set $\rho(H_{\lambda, \alpha, \theta})$ that can be checked directly. 
Finally, in Section~\ref{sec:approxSpec}, we prove results for the approximation of the spectrum of aperiodic Schr\"odinger operators with Sturmian potentials, both one- and two-sided, via their periodic counterparts.

For the remainder of this section, let us always assume the hypothesis below, which collects the notations we use.

\begin{hypo}\label{hypo:APeriodicH}
    Let $ p \in [1,\infty]$, $\alpha \in [0,1]$ be irrational, and let $\alpha_m = \frac{p_m}{q_m}$ denote the $m$-th rational approximant to $\alpha$ from Section~\ref{sec:rationalApproximations}.
    In addition, let $\theta \in [0,1)$, $\lambda\in\RR$
    and $H_{\lambda,\alpha,\theta} : \ell^p(\ZZ) \rightarrow \ell^p(\ZZ)$ be
    the aperiodic Schrödinger operator
    \begin{equation*}\label{eq:SchrodingerOp2}
    (H_{\lambda,\alpha,\theta} x)_n =  x_{n+1} +   x_{n-1} + \lambda v_{\alpha,\theta}(n) x_n \, , \quad n \in \ZZ \,,
    \end{equation*}
    with Sturmian potential
    \begin{equation*}\label{eq:defintion_words2}
            v_{\alpha,\theta} (n) = \chi_{[1-\alpha,1)}(n\alpha+\theta \bmod 1)\,.
    \end{equation*}
    Moreover, let $\widetilde{H}_{\lambda,\alpha,\theta} : \ell^p(\ZZ) \rightarrow \ell^p(\ZZ)$
    denote the aperiodic Schrödinger operator given by
    \begin{equation*}
            (\widetilde{H}_{\lambda,\alpha,\theta} x)_n =  x_{n+1} +   x_{n-1} + \lambda \widetilde{v}_{\alpha,\theta}(n) x_n \, , \quad n \in \ZZ \,,
    \end{equation*}
    where
    \begin{equation*}
            \widetilde{v}_{\alpha,\theta} (n) = \chi_{(1-\alpha,1]}(n\alpha+\theta \bmod 1) \,.
    \end{equation*}
  The sequence $(\alpha_m)$ of $m$-th rational approximants to $\alpha$ gives rise to sequences of periodic Schrödinger operators $(H_{\lambda, \alpha_m, \theta})$ and $(\widetilde{H}_{\lambda, \alpha_m, \theta})$ which are defined as for irrational $\alpha$ above. 
\end{hypo}

\subsection{The (Operator) Spectrum of Aperiodic Schr\"odinger Operators}\label{sec:limOps}

The results of Section~\ref{sec:limitWords} can be directly elevated to the level of operators, as the following result shows.
As a consequence of Proposition~\ref{prop:subshift}, we get a characterisation for the limit operators of $H_{\lambda, \alpha,\theta}$.

\begin{corollary}\label{cor:limOpsCharacterisation}
  Assume Hypothesis~\ref{hypo:words}. Then
    \begin{equation*}
  \Lla\, \coloneqq \{ H_{\lambda,\alpha,\psi}, \widetilde{H}_{\lambda,\alpha,\psi} : \psi\in[0,1) \} = \Lim(\widetilde{H}_{\lambda, \alpha,\theta}) =  \Lim(H_{\lambda, \alpha, \theta}) 
    \end{equation*}
    and
    \begin{equation*}
  \Llam\,\coloneqq \{ H_{\lambda,\alpha_m,\psi}, \widetilde{H}_{\lambda,\alpha_m,\psi} : \psi\in[0,1) \} = \Lim(\widetilde{H}_{\lambda, \alpha_m,\theta}) =  \Lim(H_{\lambda, \alpha_m, \theta}) 
  .
    \end{equation*}
  In particular, all elements of $\Lla$,
  i.e.\ all aperiodic Schrödinger operators, are self-similar.
  Moreover, the set $\Llam$ is finite. 
\end{corollary}

Note, in particular, that $\Lim(H_{\lambda, \alpha, \theta})$ and $\Lim(\widetilde{H}_{\lambda, \alpha,\theta})$ are independent of $\theta$.
Corollary \ref{cor:limOpsCharacterisation} says in short that
\begin{equation}\label{eq:minimalfamily}
\Lim(H)=\Lla\quad\text{for all}\quad H\in\Lla.
\end{equation}
In the language of dynamical systems, this fact corresponds to so-called \emph{minimality}.
It has far-reaching consequences, cf.~\cite{Lindner.2021}. Here is a first one:

\begin{prop}\label{prop:plusLimOpsEqual}
Assume Hypothesis~\ref{hypo:words}.
For all $H\in\Lla$, we have
\[
\Lla=\Lim(H)=\Lim_-(H)=\Lim_+(H)\,.
\]
\end{prop}
\begin{proof}
    Let $H\in\Lla$ and $A\in\Lim_+(H)$. Then, by~\cite[Lemma~3.3a)]{Wilde-Lindner.2008},
    $\Lim(A)\subset\Lim_+(H)$, so that by Corollary~\ref{cor:limOpsCharacterisation}
\[
\Lla=\Lim(A)\subset\Lim_+(H)\subset\Lim(H)=\Lla\,,
\]
whence both inclusions must be equalities. In particular, $\Lim_+(H)=\Lim(H)$.
The proof of $\Lim_-(H)=\Lim(H)$ is the same.
\end{proof}
Note that it is straightforward to prove that the statement of Proposition~\ref{prop:plusLimOpsEqual} also holds for all $H\in\Llam$.
We continue with some further consequences of \eqref{eq:minimalfamily}:

\begin{prop}\label{prop:lowerNormsEqual}
Assume Hypothesis~\ref{hypo:words}, and let $L \in \{\Lla,\, \Llam\}$. 
For all $A, B \in L$, the following statements hold:
    \begin{enumerate}[label=(\alph*)]
        \item $A \in \Lim(B)$ and $B \in \Lim(A)$.
        \item $\|A \| = \|B\|$ and $\nu(A) = \nu(B)$
        \item $\sigma(A) = \sigma(B)$. In particular, $A$ and $B$ are invertible at the same time and $B^{-1} \in \Lim(A^{-1})$ and $A^{-1} \in \Lim(B^{-1})$.
    \end{enumerate}

\end{prop}
\begin{proof}
$B\in\Lim(A)$ implies $\|B\|\le\|A\|$, see \cite[Proposition~1.a]{Rabinovich.1998}, $\nu(B)\ge\nu(A)$, see \cite[Proposition~3.9]{Lindner.2016} and $\sigma(B)\subset\sigma_{\ess}(A)\subset \sigma(A)$ by Lemma~\ref{prop:essSpecEqualsSpecnew}.
Since $B\in\Lim(A)$ for all $A,B\in\Lla$, see~\eqref{eq:minimalfamily}, the assertion follows.
\end{proof}

Let us turn our attention to the one-sided compressions of elements in~$\Lla$.
The next result, which is a consequence of the symmetry property~\eqref{eq:mirroredLimWord} of Sturmian words, 
will allow us to restrict ourselves to studying only one-sided compressions in $\ell^p(\ZZ_+)$.
A similar result also holds for periodic Schrödinger operators, cf.~\cite[Lemma~3.7]{Gabel.2021b}.

\begin{corollary}\label{cor:specPlusMinusEqual}
Assume Hypothesis~\ref{hypo:APeriodicH}. 
Then
   \begin{align*}
     \sigma\big( (H_{\lambda, \alpha, \theta})_+ \big)
     =
     \sigma\big( (\widetilde{H}_{\lambda, \alpha,-\alpha -\theta\bmod 1})_- \big)\,.
   \end{align*}
Moreover,
\begin{equation}\label{eq:limOpsUnion}
\bigcup_{A\in\Lla} \sigma(A_+) = \bigcup_{B\in\Lla} \sigma(B_-)\,.
\end{equation}
\end{corollary}

\begin{proof}
  Due to~\eqref{eq:mirroredLimWord}, we get $(H_{\lambda, \alpha, \theta})_+ = \Phi^{-1} (\widetilde{H}_{\lambda, \alpha,-\alpha -\theta\bmod 1})_- \Phi$ 
  with the canonical isomorphism $\Phi\colon \ell^p(\ZZ_+) \to \ell^p(\ZZ_-)$.
  This immediately yields the first equality.
  Now, equality~\eqref{eq:limOpsUnion} follows from Corollary~\ref{cor:limOpsCharacterisation}.
\end{proof}

For two-sided infinite operators $H_{\lambda, \alpha, \theta}$, we know from Proposition~\ref{prop:lowerNormsEqual} that their spectrum is independent of $\theta$. 
For the corresponding one-sided compressions, we obtain a similar result for the essential spectra.

\begin{proposition}\label{prop:essentialSpecOneSided}
Assume Hypothesis~\ref{hypo:APeriodicH} and $A,B\in\Lla$ or $A,B\in\Llam$. 
 Then,
\begin{equation*}
\sigma_\ess(A_+)= \sigma_\ess(A)=\sigma_\ess(A_-)=\sigma_\ess(B)\,.
\end{equation*}
Further,
\begin{equation*} 
\sigma(A)
    \subset 
    \sigma(A_+)
\quad \textnormal{and} \quad
\nu(A_+) \le \nu(A)\,.  
\end{equation*}
\end{proposition}
\begin{proof}
From Proposition~\ref{prop:plusLimOpsEqual} we know that 
\begin{equation*}
\Lim(A_+)= \Lim(A)=\Lim(A_-)=\Lim(H_{\lambda,\alpha, 0}) \,.
\end{equation*}
The equality of the essential spectra then follows from Proposition~\ref{prop:essSpecEqualsSpecnew}.
The rest follows from the fact that $A\in\Lim A_+$ and \cite[Proposition~3.9]{Lindner.2016}.
\end{proof}

\begin{remark}
As Example~\ref{ex:FibHamFinal} will show, the inclusion $\sigma(H_{\lambda, \alpha, \theta})
\subset 
\sigma((H_{\lambda, \alpha, \theta})_+)$ is strict. 
For examples of periodic Schrödinger operators $H$ with $\sigma(H) \subsetneq \sigma(H_+)$ see~\cite[Example~4.2, Example~4.4]{Gabel.2021b}.
\end{remark}

\subsection{Trace Condition for the Spectrum}\label{sec:traceCondition}
In this section, we characterise the invertibility of the limit operators of $H_{\lambda,\alpha,\theta}$ via traces of products of transfer matrices. 
A similar condition for the spectrum was first given by S\"ut\H{o}~\cite{Suto.1987} for the Fibonacci Hamiltonian and later generalised by Bellissard~\cite{Bellissard.1989} to Schrödinger operators with Sturmian potential. 
	
We start with the basic definitions and notation.
For an energy $\energy \in \RR$ and a vector $x \in \ker(H_{\lambda, \alpha, \theta} - E)$, we study generalised eigenfunctions, that is, solutions of 
\begin{equation}\label{eq:EVProb}
  0 = ((H_{\lambda, \alpha, \theta} - \energy)x)_n = x_{n + 1} + x_{n - 1} + (\lambda v_{\alpha, \theta} (n) - \energy) x_n\, , \quad n \in \ZZ \,.
\end{equation}
The scalar three-term recurrence in~\eqref{eq:EVProb} can be written as the vector-valued two-term recursion
\begin{equation}\label{eq:TF-Solution}
			\begin{pmatrix}
				 x_{n+1} \\
					x_n
			\end{pmatrix}	
			=
		\begin{pmatrix}
	 		 \energy-\lambda v_{\alpha,\theta}(n)  & -1 \\
   			  1 & 0
		\end{pmatrix}
		\begin{pmatrix}
			 x_n \\
			x_{n-1}
		\end{pmatrix}\,,
		\quad n \in \ZZ.
\end{equation}
The $2 \times 2$-matrix in~\eqref{eq:TF-Solution} is called \emph{transfer matrix} and here denoted 
 by
	\begin{equation}\label{eq:TF-Definition}
		T_{\lambda,\alpha,\theta}(n,\energy)  \coloneqq 
	\begin{pmatrix}
 		\energy-\lambda v_{\alpha,\theta}(n)   & -1 \\
 	     1 & 0
 	\end{pmatrix}\,,\quad n \in \ZZ\,.
	\end{equation}
          
The transfer matrices lie in the \emph{symplectic group}, which means that they satisfy
the relation~$T_{\lambda,\alpha,\theta}(n,\energy)^T \left(\begin{smallmatrix}
	0 & -1 \\ 1 & 0
  \end{smallmatrix}\right)
 T_{\lambda,\alpha,\theta}(n,\energy) 
  = \left(\begin{smallmatrix}
	0 & -1 \\ 1 & 0
  \end{smallmatrix}\right)$.
Especially, one knows that eigenvalues of a symplectic matrix always come in pairs $\{\lambda_i^{}, \lambda_i^{-1} \}$ and, therefore, their determinant is $1$.
Hence, the sequence $(x_n)_{n\in\ZZ}$ is determined by two consecutive values, say $x_0$ and $x_1$, and all other values $x_n$, $n > 0$, are determined by applying the corresponding \emph{monodromy matrix}
\begin{equation}\label{eq:TF-M}
M_{\lambda,\alpha,\theta}(n,\energy) \coloneqq  T_{\lambda,\alpha,\theta}(n,\energy) \cdots T_{\lambda,\alpha,\theta}(1,\energy)
\end{equation}
to the vector $(x_1,x_0)^\top$. 
Values $x_n$, $n < 0$, are calculated via $M_{\lambda, \alpha, \theta}(n, \energy)^{-1}$.

The following proposition resembles a well-known result about the resolvent of periodic Schrödinger operators, which can be described by a trace condition of the monodromy matrix, cf., e.g.,~\cite{Puelz.2014}.

\begin{proposition}\label{prop:trace_condition1}
  Assume Hypothesis~\ref{hypo:APeriodicH}, let $\energy \in \RR$, and consider $\alpha_m = \frac{p_m}{q_m}$. 
  Then $\energy \in \rho(H_{\lambda, \alpha_m, \theta})$ if and only if the so-called \emph{trace condition}
  \begin{equation*}
  \abs{\tr(M_{\lambda, \alpha_m,\theta}(q_m,E))} > 2
  \end{equation*}
    holds.
\end{proposition}
	
The following proposition gives a description of the resolvent set of aperiodic Schrödinger operators in terms of traces of $M_{\lambda,\alpha,\theta}(n,\energy)$. 
Note that, in contrast to Proposition~\ref{prop:trace_condition1}, it is not sufficient to consider only one matrix $M_{\lambda,\alpha,\theta}(n,\energy)$.

\begin{proposition}[{\cite[Proposition 4]{Bellissard.1989}}]\label{prop:SpurKriterium}
  Assume Hypothesis~\ref{hypo:APeriodicH}, and let $E\in\RR$. 
  Then $\energy \in\rho(H_{\lambda,\alpha,\theta})$ if and only if the \emph{trace condition}
  \begin{align}\label{eq:traceCondAperiodicH}
    \exists m\geq 0 : \,
    \left|\tr\big(M_{\lambda,\alpha,0}(q_m, \energy) \big)\right| 
    > 2 
    \; \text{and} \;
    \left| \tr \big( M_{\lambda,\alpha,0}(q_{m+1}, \energy) \big)\right| 
    > 2
  \end{align}
  holds, where $q_m\in \NN$ is given by the $m$-th rational approximant $\frac{p_m}{q_m}$ to $\alpha$ for $m\in \NN$.
In this case, 
  \begin{align*}
    \left| \tr \big( M_{\lambda,\alpha,0}(q_k, \energy) \big)\right| > 2
  \end{align*}
for all $k\ge m$.
\end{proposition}

Note that the trace condition~\eqref{eq:traceCondAperiodicH} in Proposition~\ref{prop:SpurKriterium} only considers the monodromy matrices with $\theta = 0$ while the conclusion about the resolvent set is for general $\theta$. 
This reflects the fact that $\sigma(H_{\lambda,\alpha,\theta})$ is independent of $\theta\in [0,1)$ as shown Proposition~\ref{prop:lowerNormsEqual}. 

\begin{remark}\label{rem:traces1}
  In the trace condition~\eqref{eq:traceCondAperiodicH} in Proposition~\ref{prop:SpurKriterium}, we can also use the monodromy matrices of the periodic Schrödinger operators $H_{\lambda,\alpha_m,0}$ and $H_{\lambda,\alpha_{m+1},0}$ for calculating the trace because of the relation
        \begin{equation}\label{eq:tracePeriodicAperiodic}
        \tr( M_{\lambda,\alpha,0}(q_m, \energy )) = \tr( M_{\lambda,\alpha_m,0}(q_m, \energy )) \,.
        \end{equation}
  Let us verify~\eqref{eq:tracePeriodicAperiodic} for $m$ even and $m$ odd separately.

  For $m=2$ this follows directly from Lemma~\ref{lem:potential_aperiodic1}.
  For $m \geq 4$ even, $\trans_{\lambda, \alpha_m, 0}(k, \energy) = \trans_{\lambda, \alpha, 0}(k, \energy)$ for all $k \in \{1,\dots,q_m\}$ as a consequence of Lemma~\ref{lem:potential_aperiodic1}.
        The equality of the transfer matrices then directly implies relation~\eqref{eq:tracePeriodicAperiodic}.

        For $m \geq 3$ odd, we use the translation of the potential given in 
	Remark~\ref{rem:shift_invariant} with~$k = q_m$ and $n=m - 1$ which allows the calculation
	\begin{align*}
		M_{\lambda,\alpha,0}(q_m,\energy )
		&= 
		\trans_{\lambda, \alpha, 0}(q_m, \energy) \cdots \trans_{\lambda, \alpha, 0}(1, \energy)  \\
		&= 
		\trans_{\lambda, \alpha, 0}(q_m - q_{m-1}, \energy) \cdots \trans_{\lambda, \alpha, 0}(1 - q_{m-1}, \energy)    \\
		&= 
		\trans_{\lambda, \alpha_m, 0}(q_m - q_{m-1}, \energy) \cdots \trans_{\lambda, \alpha_m, 0}(1 - q_{m-1}, \energy) \,.
	\end{align*}
	The last equality holds again by Lemma~\ref{lem:potential_aperiodic1}. 
        Now, taking the trace on both sides and using the periodicity of $v_{\alpha_m}$ while reordering the product, we arrive also in this case at the claimed relation~\eqref{eq:tracePeriodicAperiodic}. 
\end{remark}

In~\cite{Bellissard.1989}, Proposition~\ref{prop:SpurKriterium} was only shown for positive coupling constant~$\lambda$. 
It can, however, be extended to $\lambda\in\RR$, where $\lambda = 0$ is just the periodic case.
We will postpone this extension to Remark~\ref{rem:traceConditionNegative}

In the following corollary, we connect Proposition~\ref{prop:SpurKriterium} to the invertibility problem of all limit operators of a given aperiodic Schrödinger operator.
\begin{corollary}\label{cor:LimOpsInvIffTrace}
Assume Hypothesis~\ref{hypo:APeriodicH}.
All limit operators of $H_{\lambda, \alpha, \theta}$ are invertible if and only if
  \begin{equation}\label{eq:traceCondition}
    \exists m\geq 0 : \;
    \left|\tr\big(M_{\lambda,\alpha,0}(q_m, 0) \big)\right| 
    > 2 
    \;\text{and}\;
    \left| \tr \big( M_{\lambda,\alpha,0}(q_{m+1}, 0) \big)\right| 
    > 2
\,.
  \end{equation}
\end{corollary}

\begin{proof}
Let $A\in \operatorname{Lim}(H_{\lambda,\alpha,\theta})$. 
Due to Corollary~\ref{cor:limOpsCharacterisation}, there exists $\psi\in [0,1)$ such that either $A=H_{\lambda,\alpha,\psi}$ or $A=\widetilde{H}_{\lambda,\alpha,\psi}$. 
By Proposition~\ref{prop:lowerNormsEqual}, we have $\sigma(H_{\lambda,\alpha,\theta}) = \sigma(H_{\lambda,\alpha,\psi}) = \sigma(\widetilde{H}_{\lambda,\alpha,\psi})$. 
  Consequently, $A$ is invertible if and only if $H_{\lambda, \alpha, \theta}$ is invertible. Since by Proposition~\ref{prop:kurbatov} the spectrum is independent of $p\in [1,\infty]$, the claim then follows by Proposition~\ref{prop:SpurKriterium}.
\end{proof}

We close this section by verifying the trace condition~\eqref{eq:traceCondition} for the Fibonacci Hamiltonian with energy $E=0$.

\begin{example}\label{ex:FibHam}
Assume Hypothesis~\ref{hypo:APeriodicH}.
Let $\alpha = (\sqrt{5}-1)/2$ and $\lambda \in \RR$. 
 Then,
  \begin{align*}
    \tr(M_{\lambda,\alpha,0}(q_6,0)) &= \lambda  ( 18 \lambda^2-8 )\,,  \\
    \tr(M_{\lambda,\alpha,0}(q_7,0)) &=\lambda(-108 \lambda^4 +84 \lambda^2-13) \,.
  \end{align*}
For $\lambda > \frac{1}{6} (3 + \sqrt{3}) \eqqcolon \lambda_0 \approx 0.788675$, we have
  \begin{equation*}
    \tr\big(M_{\lambda,\alpha,0}(q_6,0)\big) > 2 
    \quad\text{and}\quad
    \tr\big(M_{\lambda,\alpha,0}(q_7,0)\big) < -2 \,. 
  \end{equation*}
Hence, by Corollary~\ref{cor:LimOpsInvIffTrace}, for $\lambda > \lambda_0$, all limit operators of the Fibonacci Hamiltonian $H_{\lambda,\alpha,0}$ are invertible. 
The bound~$\lambda_0$ may be further improved  by considering higher orders of approximation. 
Considering the polynomial~$\tr(M_{\lambda,\alpha,0}(q_8,0))$, however, does not improve the result despite being algebraically solvable as well.
  Figure~\ref{fig:FibHamLambda} visualises the spectra of the periodic approximations for $m = 5,6,7,8$ with special emphasis on the previously derived bound $\lambda_0$, see also~\cite[Section~7.1]{Damanik.2015}.

\begin{figure}[htbp]
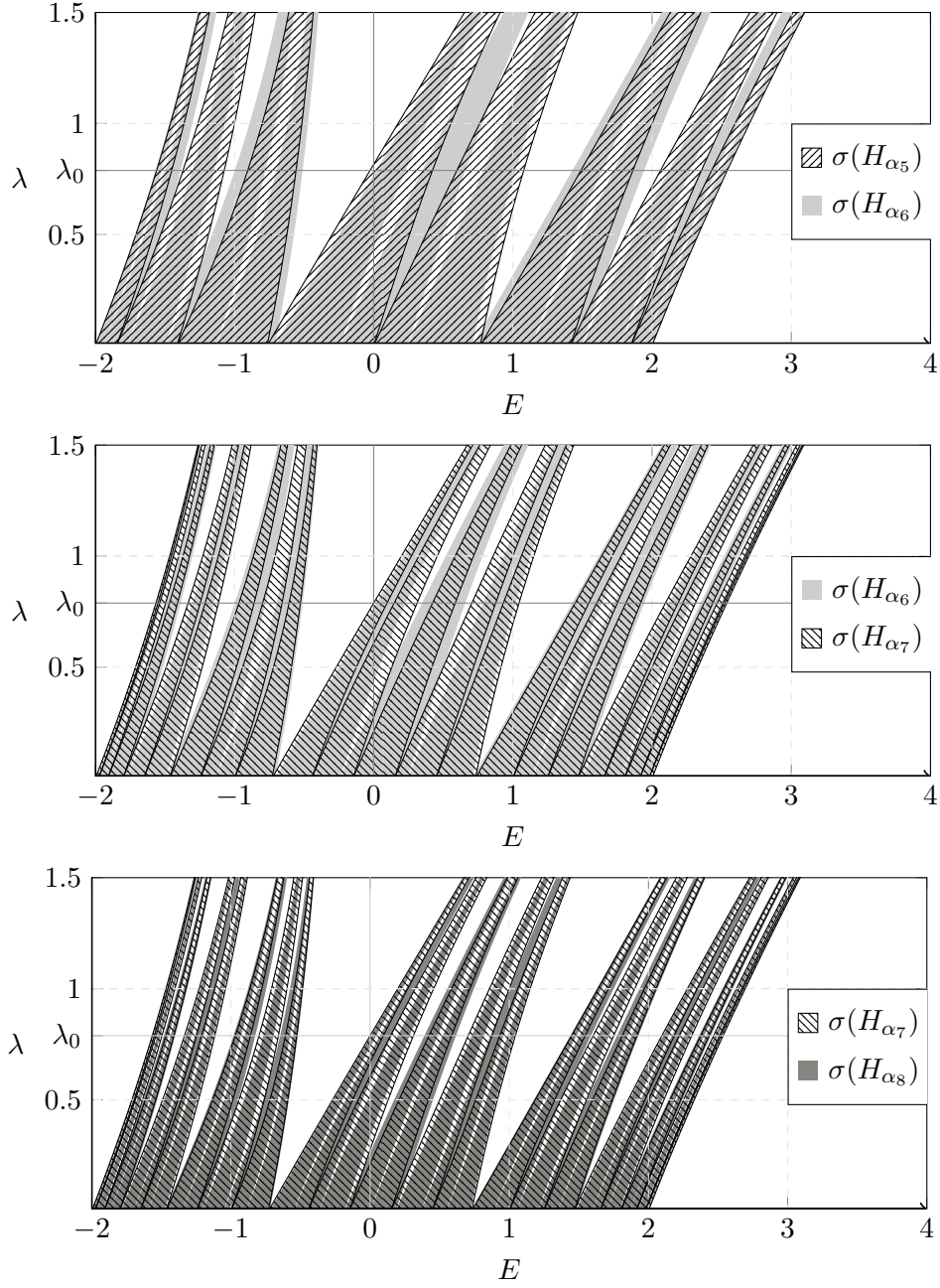
 
    \def\drawqfve{\draw[ pattern=north east lines]}
\def\drawqsx{\fill[bandSpecCol!40]}
\def\drawqsvn{\draw[pattern=north west lines]}
\def\drawqeight{\fill[bandSpecCol]}

     \caption{Spectra of the periodic approximations~$H_{\lambda,\alpha_m,0}$ for~$m = 5,6,7,8$ of the scaled Fibonacci Hamiltonian~$H_{\lambda,\alpha,0}$, where $\alpha=(\sqrt{5} - 1)/2$ and $\lambda >0$. 
    If $\lambda > \lambda_0$, then $0 \not \in \sigma(H_{\lambda, \alpha, 0})$ as $\sigma(H_{\lambda,\alpha,0}) \subset \sigma(H_{\lambda, \alpha_m,0}) \cup \sigma(H_{\lambda,\alpha_{m + 1},0})$ by Corollary~\ref{prop:SpurKriterium} and all limit operators of $H_{\lambda, \alpha, 0}$ are invertible by Corollary~\ref{cor:LimOpsInvIffTrace}.}\label{fig:FibHamLambda}
\end{figure}
\end{example}

\subsection{Approximation of the Spectrum} \label{sec:approxSpec}
In this section, we introduce the approximation theory for both one- and two-sided aperiodic Schrödinger operators.
As described in Hypothesis~\ref{hypo:APeriodicH}, every irrational number $\alpha$ with rational approximants $(\alpha_m)$ gives rise to \emph{aperiodic} Schrödinger operators $H_{\lambda, \alpha, \theta}$, $\widetilde{H}_{\lambda, \alpha, \theta}$ and sequences $(H_{\lambda, \alpha_m, \theta})$, $(\widetilde{H}_{\lambda, \alpha_m, \theta})$ of \emph{periodic} Schrödinger operators.
The following proposition, which originally appeared in~\cite{Suto.1987}, shows in which sense the aforementioned approximation sequences can be considered a \emph{periodic approximation} of aperiodic Schrödinger operators.  

\begin{prop}[{\cite[Proposition~4]{Suto.1987}}]\label{prop:rationalApproximation}
Assume Hypothesis~\ref{hypo:APeriodicH}.
Then, 
\begin{align*}
  H_{\lambda,\alpha,0}\, x = \lim_{m\to\infty} H_{\lambda,\alpha_m,0} \, x \qquad\text{and}\qquad
  \widetilde{H}_{\lambda,\alpha,0}\, x = \lim_{m\to\infty} \widetilde{H}_{\lambda,\alpha_m,0} \, x
\end{align*}
for all $x \in \cfin(\ZZ)$.
Furthermore, if $p < \infty$,
then 
\begin{align}\label{eq:strongConvergencePeriodicApprox}
  H_{\lambda,\alpha,0} = \lim_{m\to\infty} H_{\lambda,\alpha_m,0} \qquad\text{and}\qquad
  \widetilde{H}_{\lambda,\alpha,0} = \lim_{m\to\infty} \widetilde{H}_{\lambda,\alpha_m,0}
\end{align}
in the strong operator topology.
\end{prop}

\begin{proof}
  We only prove the claim for the sequence $(H_{\lambda, \alpha_m, 0})$.
  The proof for $(\widetilde{H}_{\lambda, \alpha_m, 0})$ is completely analogous.

  Recall that the rational number $\alpha_m$ induces the periodic two-sided word $v_{\alpha_m}$.
  Furthermore, by Lemma~\ref{lem:potential_aperiodic1}, we know that $v_{\alpha}(k) = v_{\alpha_m}(k)$ for all $k \in \ZZ$ with $|k| \leq q_m - 2$ and $m \geq 2$.
This shows the first part.
The second part then follows from the fact that $\cfin(\ZZ)$ is dense in $\ell^p(\ZZ)$ for $p<\infty$.
\end{proof}

It is easy to see that the assertion of Proposition~\ref{prop:rationalApproximation} also holds mutatis mutandis for the one-sided compressions $( H_{\lambda,\alpha,0})_+$ and  $( H_{\lambda,\alpha,0})_-$ and, for a certain choice of subsequences $(\alpha_{m_k})$, also for $H_{\lambda,\alpha,\theta}$ and $(H_{\lambda,\alpha_{m_k},\theta})$, $\theta\in [0,1)$, see Lemma~\ref{lem:potential_aperiodic_theta}.

In the following, we analyse how the approximation result of Proposition~\ref{prop:rationalApproximation} translates into convergence of the respective spectra.
We start with a general approximation result.
When stating results concerning one-sided operators, recall that, due to Corollary~\ref{cor:specPlusMinusEqual}, we only need to treat compressions $A_+$.

Similar to Section~\ref{sec:limitWords}, for a Schrödinger operator $A$ with diagonal $v$, let $\mathcal{W}_D(A)$, 
$D \in \NN$, denote the list of all subwords of the diagonal of $A$ that have length $D$, 
i.e.~$\mathcal{W}_D(A) = \mathcal{W}_D(v)$. 
Furthermore, let $\mathcal{W}(A) \coloneqq \bigcup_{D \in \NN} \mathcal{W}_D(A)$ denote the list of all finite subwords of the diagonal of $A$.

\begin{lemma}\label{lem:subwordsAndLowerNorms}
  Let $A$, $A_m$, $m \in \NN$, and $B$ be Schrödinger operators on $\ell^p(\ZZ)$, $p \in [1, \infty]$.
  \begin{enumerate}[label=(\alph*)]
      \item\label{it:sameSubwordsApprox} 
        If $\mathcal{W}_D(A) \subset \mathcal{W}_D(B)$ for some $D \in \NN$, then $\nu_D(A) \geq \nu_D(B)$.
        In particular, if $\mathcal{W}(A) = \mathcal{W}(B)$, then $\nu(A)= \nu(B)$.
      \item\label{it:decomposition}
        If $\mathcal{W}_D(A) \subset \mathcal{W}_D(A_+)$, then
        \begin{equation}\label{eq:lowerNormPlus}
          \nu_D(A_+)= \min\Big\{ \nu_{0..D}(A_+),\;\;  \nu_D(A) \Big\} \,,
        \end{equation}
      where $\nu_{k..l}(A) \coloneqq \inf \big\{ \|Ax \| : \supp(x) \subset \{k,\dots,l\}, \|x\| = 1 \big\}$.
      \item\label{it:limit}
          Assume that $(A_m)$ is uniformly bounded and, for each $D \in \NN$, we have $\mathcal{W}_D(A_m) = \mathcal{W}_D(A)$ eventually.
          Then
          \begin{align*}
            \lim_{m \to \infty} \nu(A_m) = \nu(A)\,.
          \end{align*}
      \item\label{it:limitOneSided}
        If, in addition to~\ref{it:limit}, $(A_m)_+x \to A_+x$ for all $x\in \cfin(\ZZ_+)$, we have 
        \begin{align*}
          \lim_{m \to \infty} \nu((A_m)_+) = \nu(A_+)\,.
        \end{align*}
  \end{enumerate}
\end{lemma}

\begin{proof}
  \ref{it:sameSubwordsApprox}: Let $x \in \ell^p(\ZZ)$ with $\|x\|=1$ and $\diam( \supp(x)) \leq D$. 
  As $\mathcal{W}_D(A) \subset \mathcal{W}_D(B)$, there exists $k \in \ZZ$ such that $\|A x \| = \|B S^k x \|$. This implies $\nu_D(A) \geq \nu_D(B)$.
    The second part of the claim follows by a symmetry argument and the fact that $\nu(A) = \inf_{D \in \NN} \nu_D(A)$, see Lemma~\ref{lem:lowerNorm}.

    \ref{it:decomposition}:
    As the approximate lower norm~$\nu_D$ considers only sequences with finite support, the following decomposition is evident:
   \begin{align*}
     \nu_D(A_+) 
     = \min \Big\{ \nu_{0..D}(A_+),\;\; \inf_{d \in \NN} \nu_{d..d+D}(A_+) \Big\} \,.
   \end{align*}  
    In order to verify identity~\eqref{eq:lowerNormPlus}, we need to show that
    \begin{align}\label{eq:lowerNormPositive}
      \nu_D(A) = \inf_{d \in \NN} \nu_{d..d+D}(A_+)\,.
    \end{align}
    To this end, let $x \in \cfin(\ZZ)$ with $\diam(\supp(x)) \leq D$. 
    As by assumption $\mathcal{W}_D(A) \subset \mathcal{W}_D(A_+)$, there exists $k \in \NN$ such that 
    \begin{align*}
      \supp(S^k x) \subset \NN \quad\text{and}\quad  A x = A S^k x = A_+ S^k x \,.
    \end{align*}
    This gives that 
    \begin{equation*}
      \nu_D(A) \geq \inf_{d \in \NN} \nu_{d..d + D}(A_+)\,.
    \end{equation*}
    The estimate $\nu_D(A) \leq \nu_{d .. d + D}(A_+)$ trivially holds for all $d \in \NN$.
    This proves~\eqref{eq:lowerNormPositive}.

    \ref{it:limit}: 
    Let $\varepsilon > 0$ and take $D$ from Lemma~\ref{lem:lowerNorm}. 
    By part~\ref{it:sameSubwordsApprox},
    \begin{align}\label{eq:nuEstimatesCombined}
      \nu(A) \leq \nu_D(A_m) \leq \nu(A_m) + \varepsilon 
      \quad\text{and}\quad
      \nu(A_m) \leq\nu_D(A) \leq \nu(A) + \varepsilon
    \end{align}
    eventually.
    As $\varepsilon > 0$ was arbitrary, the first part of~\eqref{eq:nuEstimatesCombined} gives
    \begin{align*}
        \nu(A) \leq \liminf_{m \to \infty} \nu(A_m)
    \end{align*}
    while the second part of~\eqref{eq:nuEstimatesCombined} gives
    \begin{align}\label{eq:limsupAm}
        \limsup_{m \to \infty} \nu(A_m) \leq \nu(A)\,. 
    \end{align}

    \ref{it:limitOneSided}:
    From $(A_m)_+x \to A_+x$ for $x\in\cfin(\ZZ_+)$, we directly get that
    \begin{align*}
      \nu_{0..D}((A_m)_+) = \nu_{0..D}(A_+)
    \end{align*}
    eventually. 
    Furthermore, by part~\ref{it:sameSubwordsApprox},
    \begin{align*}
      \nu_D(A_m) = \nu_{D}(A)
    \end{align*}
    eventually.
    Now, part~\ref{it:decomposition} gives
    \begin{align*}
      \nu_D((A_m)_+) 
      &= \min \Big\{ \nu_{0..D}((A_m)_+),\;\; \nu_D(A_m) \Big\}  \\
      &= \min \Big\{ \nu_{0..D}(A_+),\;\; \nu_D(A) \Big\} 
      = \nu_D(A_+) 
    \end{align*}
    eventually.
    Arguing as in the proof of part~\ref{it:limit}, the claim follows.
\end{proof}

\begin{remark} 
  By inspecting the proof of Lemma~\ref{lem:subwordsAndLowerNorms}\ref{it:limit}, one gets, for the case in which, for all $D \in \NN$, only one of the assumptions 
  \begin{equation*}
    \mathcal{W}_D(A_m) \subset \mathcal{W}_D(A) \quad\text{or}\quad \mathcal{W}_D(A) \subset \mathcal{W}_D(A_m)
  \end{equation*}
  holds, eventually the estimate
  \begin{equation*}
    \nu(A) \leq \liminf_{m \to \infty} \nu(A_m) \quad\text{or}\quad \limsup_{m \to \infty} \nu(A_m) \leq \nu(A)\,,
  \end{equation*}
  respectively.
  For more details on this observation see \cite{Gabel.2022}.
\end{remark}

In the case of Schrödinger operators $H_{\lambda, \alpha, \theta}$ and $H_{\lambda, \alpha_m \theta}$, Lemma~\ref{lem:subwordsAndLowerNorms} gives the following result.

\begin{prop}\label{prop:specApprox}
Assume Hypothesis~\ref{hypo:APeriodicH}.
Set $A= H_{\lambda, \alpha,\theta}$ and $A_m =H_{\lambda, \alpha_m,\theta}$, and let $E\in\RR$. 
Then,
\begin{enumerate}[label=(\roman*)]
 \item\label{it:LNapproxPlus} $ \nu (A-E)  = \lim_{m\rightarrow\infty}\nu (A_m-E)$,
  \item\label{it:Stableapprox}  $(A_m-E)_{m\in\NN}$ is stable if and only if $A-E$ is invertible. In which case,
    $\lim_{n\rightarrow\infty} \|(A_m-E)^{-1}\| = \| (A-E)^{-1}\|$. 
  \item \label{it:Specapprox}$ \sigma (A_m)\to \sigma (A)$ as $m\to\infty$ in the Hausdorff metric.
\end{enumerate}
\end{prop}
\begin{proof}
Part~\ref{it:LNapproxPlus} follows directly from Lemma~\ref{lem:sameSubwordsApprox} and Lemma~\ref{lem:subwordsAndLowerNorms}\ref{it:limit}.
Part~\ref{it:Stableapprox} follows from part~\ref{it:LNapproxPlus} by using Lemma~\ref{lem:approxMerged}.
Now part~\ref{it:Specapprox} is standard, see also~\cite{Lindner.2023}.
\end{proof}
Recall that the \emph{Hausdorff  distance} of two compact sets $S,T\subset\CC$ is defined by
\[
d_{\rm H}(S,T)\ :=\ \min\left( \inf_{s\in S} \dist(s,T)\,,\, \inf_{t\in T}\dist(t,S)\right)
\]
and that $d_{\rm H}$ is a metric on the set of all compact subsets of $\CC$.
\begin{remark}\label{rem:traceConditionNegative}
In~\cite{Bellissard.1989}, Proposition~\ref{prop:SpurKriterium} was only shown for positive coupling constant~$\lambda$. 
The extension to $\lambda \leq 0$, where $\lambda = 0$ is just the periodic case, works as follows.
Let $\lambda < 0$, and observe that $T_{-\lambda,\alpha,\theta}(k,-\energy )=-T_{\lambda,\alpha,\theta}(k,\energy )^T$.
As a consequence of the simplecticity of the monodromy matrices, one gets the identity
  \begin{equation*}
    \left|\tr\big(M_{-\lambda,\alpha,0}(q_m,-\energy)\big)\right| 
    = 
    \left|\tr\big(M_{\lambda,\alpha,0}(q_m,\energy)\big)\right|\,.
  \end{equation*}
See~\cite[Lemma~3.7]{Gabel.2021b} for details.

Due to Remark~\ref{rem:traces1} and Proposition~\ref{prop:trace_condition1}, we have $\sigma(H_{\lambda, \alpha_m, \theta}) = - \sigma(H_{-\lambda, \alpha_m, \theta})$. 
With Proposition~\ref{prop:specApprox}, the identity of spectra carries over to the aperiodic operators. 
Thus, we also have
  \begin{equation*}
- \sigma(H_{\lambda, \alpha, \theta}) = \sigma(H_{-\lambda, \alpha, \theta})\,,
  \end{equation*}
  i.e.~the spectrum for $H_{-\lambda, \alpha, \theta}$ is already determined by $H_{-\lambda, \alpha, \theta}$ which in turn is characterised by Proposition~\ref{prop:SpurKriterium}.
\end{remark}

In order to get result similar to Proposition~\ref{prop:specApprox} for the one-sided case, we have to assure that the condition $(A_m)_+x \to A_+x$ for $x\in\cfin(\ZZ_+)$ in Lemma~\ref{lem:subwordsAndLowerNorms}\ref{it:limitOneSided} is satisfied. 
This requires a careful choice of sequences $A_m$ depending on $\theta$.
Fortunately, Lemma~\ref{lem:potential_aperiodic_theta} yields the following sequences:

\begin{table}[h!]
\begin{center}
{\def\arraystretch{1.5}\tabcolsep=1.3ex
\begin{tabular}{ r|c|c| } 
     &  $A = H_{\lambda, \alpha, \theta}$ & $A = \widetilde{H}_{\lambda, \alpha, \theta}$ \\ 
 \hline
   $\theta = 0$ & $H_{\lambda, \alpha_m, 0}$ & $\widetilde{H}_{\lambda, \alpha_m, 0}$ \\
   $\theta = -k_0 \alpha \bmod 1$ with $k_0 \in \NN$ &  $H_{\lambda, \alpha_{2m + 1}, \theta}$ & $H_{\lambda, \alpha_{2m}, \theta}$  \\
   $\theta \neq -k_0\alpha \bmod 1$ for all $k_0 \in \NN$ & $H_{\lambda, \alpha_m, \theta}$ & $H_{\lambda, \alpha_m, \theta}$  \\ 
 \hline
\end{tabular}
}
\end{center}
  \caption{Given a periodic Schrödinger operator $A$ and a condition on the offset $\theta$, the corresponding table entry gives sequence $(A_m)$ used for approximation.}
\label{table:choiceAm}
\end{table}

\begin{prop}\label{prop:specApproxPlus}
Assume Hypothesis~\ref{hypo:APeriodicH}.
Let $A$ and $A_m$ be defined as in Table~\ref{table:choiceAm} and $E\in\RR$. 
Then
\begin{enumerate}[label=(\roman*)]
 \item \label{it:LNapproxPlus} $ \nu (A_+)  = \lim_{m\rightarrow\infty}\nu ((A_m)_+)$,
  \item \label{it:StableapproxPlus} $((A_m-E)_+)_{m\in\NN}$ is stable if and only if $(A-E)_+$ is invertible,
  \item \label{it:SpecapproxPlus} $\sigma ((A_m)_+)\to \sigma (A_+)$ as $m\to\infty$ in the Hausdorff metric.
\end{enumerate}
Moreover, if $E\not\in \sigma(A)$ then
  \begin{align*}
  \left\| (A-E)^{-1}_+\right\| =\lim_{m\rightarrow\infty} \left\|(A_m-E)^{-1}_+\right\|\,. 
  \end{align*}  
\end{prop}

\begin{proof}
It follows from Lemma~\ref{lem:potential_aperiodic1}, Lemma~\ref{lem:potential_aperiodic_theta}, and Lemma~\ref{lem:sameSubwordsApprox} that the conditions in Lemma~\ref{lem:subwordsAndLowerNorms}\ref{it:limitOneSided} are satisfied.
Hence, we immediately arrive at \ref{it:LNapproxPlus}.

Part \ref{it:StableapproxPlus} and \ref{it:SpecapproxPlus} follow exactly as in Proposition~\ref{prop:specApprox}.
\end{proof}

The results of Proposition~\ref{prop:specApproxPlus} allow us to check the invertibility of all operators~$(H_{\lambda, \alpha, \theta})_+$ for $\theta\in[0,1)$. 
The following proposition shows that it is sufficient to only check this for $\theta\in \{k_0\alpha \bmod 1 : k_0\in\ZZ\}$ if a uniform boundedness condition is satisfied.

\begin{prop}\label{prop:onlyOrbit}
	Assume Hypothesis~\ref{hypo:APeriodicH}. Then	
\begin{equation}\label{eq:LNorbit}
 \inf_{B\in\Lla}\nu(B_+)  =  \min_{B\in\Lla}\nu(B_+)  = \inf_{k\in\ZZ}  \nu((H_{\lambda, \alpha, \theta_k})_+) \,.
\end{equation}
In particular, all operators $B_+$ for $B\in\Lla$ are invertible if and only if for all $k\in\ZZ$ the operators $(H_{\lambda, \alpha, \theta_k})_+$ with $\theta_k\coloneqq k\alpha \bmod 1$ are invertible and their inverses are uniformly bounded in $k$.
 In that case,
 \begin{equation}\label{eq:InvNormApprox}
\sup_{B\in\Lla}\big(\| B_+^{-1}\| \big) = \sup_{k\in\ZZ} \| (H_{\lambda, \alpha, \theta_k})_+^{-1}\| 
 \end{equation}
 \end{prop}
\begin{proof}
We  only show \eqref{eq:LNorbit}.
  For the first identity, take a sequence $(A_n)\subset \Lla$ such that $\lim_{n \to \infty} \nu((A_n)_+) =  \inf_{B\in\Lla}\nu(B_+) $.
By Lemma~\ref{lem:thetaConvergence}\ref{it:orbitTF1}, there exists $C\in\Lla$ such that for every $x\in\cfin(\ZZ_+)$ $A_n x\to Cx$ as $n\to \infty$. 
Together with Lemma~\ref{lem:subwordInRightHalf}, we find that the conditions of Lemma~\ref{lem:subwordsAndLowerNorms}\ref{it:limitOneSided} are satisfied and therefore $\nu(C_+) = \inf_{B\in\Lla}\nu(B_+)$. 
 Hence, the infimum in~\eqref{eq:LNorbit} is indeed a minimum.

Since $\{H_{\lambda, \alpha, \theta_k} : k\in\ZZ\}$ is a subset of $\Lla$, we immediately get 
\begin{equation*}
\min_{B\in\Lla}\nu(B_+) \leq \inf_{k\in\ZZ} \nu((H_{\lambda, \alpha, \theta_k})_+)\,. 
\end{equation*}

For the other inequality, take $C\in\Lla$ with $\nu(C_+)=\min_{B\in\Lla}\nu(B_+)$. 
  By Lemma~\ref{lem:thetaConvergence}\ref{it:orbitTF2}, we find a subsequence $(\theta_{k_j})$ such that $C_+x = \lim_{j \to \infty} (H_{\lambda, \alpha, \theta_{k_j}})_+x$ for all $x\in\cfin(\ZZ_+)$.
  As $(H_{\lambda, \alpha, \theta_k})$ is uniformly bounded, by Lemma~\ref{lem:subwordsAndLowerNorms}\ref{it:limitOneSided}, we find $\nu((H_{\lambda, \alpha, \theta_{k_j}})_+) \rightarrow \nu(C_+) $ which proves that  
\begin{equation*}
 \min_{B\in\Lla}\nu(B_+)   \geq \inf_{k\in\ZZ}  \nu((H_{\lambda, \alpha, \theta_k})_+) \,. \qedhere
\end{equation*}
 \end{proof}

\begin{remark}\label{rem:onlyOrbit2}
  With the same method as in the proof of Proposition~\ref{prop:onlyOrbit}, one can also check that $(H_{\lambda, \alpha, \theta})_+$ and $(\widetilde{H}_{\lambda, \alpha, \theta})_+$ are invertible for every $\theta \in [0,1)$ if and only if, for all $k\in\ZZ$, the operators $(\widetilde{H}_{\lambda, \alpha, \theta_k})_+$ are invertible and their inverses are uniformly bounded in $k$. 
 
 In particular, we have that the operators $(H_{\lambda, \alpha, \theta_k})_+$, $k\in\ZZ$, are all invertible and their inverses are uniformly bounded in $k$ if and only if the same holds for $(\widetilde{H}_{\lambda, \alpha, \theta_k})_+$.
\end{remark}

The results of Section~\ref{sec:specHltheta} have shown that both the spectrum of one- and two-sided aperiodic Schrödinger operators with Sturmian potential can be approximated via the spectra of their respective periodic approximations.
Since the periodic approximations are periodic Schrödinger operators, their spectra can be explicitly computed, as shown in the following result that is based on~\cite[Theorem~4.42(i)]{Hagger.2016} and~\cite[Proposition~3.8]{Gabel.2021b}.

\begin{proposition}\label{prop:oneSidedSpecRep}
    Let $H \colon \ell^p(\ZZ) \to \ell^p(\ZZ)$ be a periodic Schrödinger operator with period $K \geq 2$.
    Then $\CC \setminus \sigma_\ess(H_+)$ contains no bounded connected components and
    \begin{equation}\label{eq:oneSidedSpecRep}
                \sigma(H_+) = \sigma(H) \cup \Big\{E \in \sigma\Big((H_{i,j})_{i,j = 0}^{\period - 2}\Big) 
                : 
                | \det\Big(((H - EI)_{i,j})_{i,j = 0}^{\period - 1} \Big)| < 1 \Big\} \,.
    \end{equation}
    In particular, $H$ is invertible if $H_+$ is invertible.
\end{proposition}
Now we are in the position to carry over these observations to applicability conditions for the FSM.

\section{Finite Section Method for \texorpdfstring{$H_{\lambda, \alpha,\theta}$}{H} }\label{sec:fsm}

In the previous section, we have used an approximation method with periodic operators that has been properly tailored for the aperiodic Schrödinger operators. By using
the periodic approximations, we are able to approximate spectra and lower norms of $H_{\lambda,\alpha,\theta}$.
However, these approximants are still operators of infinite rank.
Therefore, when it comes to solving a linear system $H_{\lambda,\alpha,\theta} x = b$ for $x,b\in\ell^p(\ZZ)$, the FSM is better suited since it allows to pass to  finite-dimensional systems, which are usually easier to solve.

\subsection{Characterisation of Applicability}\label{sec:applicability}

Recall from Proposition~\ref{prop:FSMCondition} that the applicability of the FSM depends on the invertibility of all one-sided limit operators. 
Therefore, the results of Section~\ref{sec:limOps} about the spectrum of these operators
immediately translate to conditions for the applicability of the FSM.

\begin{prop}\label{prop:FSMiffPlus}
	Let $ p \in [1,\infty]$, $\alpha \in [0,1]$ be irrational, $\theta \in [0,1)$, and $\lambda\in\RR$. Then
the FSM for $H_{\lambda, \alpha, \theta}$ is applicable if and only if the operators $(H_{\lambda, \alpha, \psi})_+$ and  $(\widetilde{H}_{\lambda, \alpha, \psi})_+$ are invertible for every $\psi\in [0,1)$. The same equivalence holds for the applicability of the FSM for $\widetilde{H}_{\lambda, \alpha, \theta}$ .
\end{prop}
\begin{proof}
Let us the abbreviation $A\coloneqq H_{\lambda, \alpha, \theta}$.
From Proposition~\ref{prop:FSMCondition}, we know that the FSM is applicable if and only if all the following operators are invertible:

\begin{enumerate}[label=(\alph*)]
\item\label{it:opA1} The operator $A$, \qquad
\item\label{it:opB-1} the compressions $\,B_-$ for all $B\in  \Lim_+(A)$, and \qquad
\item\label{it:opC+1} $\,C_+$ for all $C\in  \Lim_-(A)$.
\end{enumerate}
Recall from Corollary~\ref{cor:limOpsCharacterisation} and Proposition~\ref{prop:plusLimOpsEqual} that the limit operators of $A$, as well as $\Lim_+(A)$ and $\Lim_-(A)$, are given by $\{ H_{\lambda,\alpha,\psi}, \widetilde{H}_{\lambda,\alpha,\psi} : \psi\in[0,1) \}$. Moreover,
from Corollary~\ref{cor:specPlusMinusEqual} and Proposition~\ref{prop:plusLimOpsEqual} we know that the invertibility of all operators in~\ref{it:opC+1} implies the invertibility of all operators in~\ref{it:opB-1}. 
According to Proposition~\ref{prop:essentialSpecOneSided}, the invertibility of $C_+$ implies the invertibility of all limit operators of $A$. 
In particular, this implies the invertibility of $A$, since $A\in\Lim(A)$ by Corollary~\ref{cor:limOpsCharacterisation}. 
Hence, the invertibility of all operators in~\ref{it:opC+1} is already equivalent to the applicability of the FSM. 
\end{proof}

Now we are able to present a collection of different conditions that guarantee the applicability of the FSM for the aperiodic Schrödinger operator.
Here, in the case that an operator $A$ is not invertible, we simply define $\|A^{-1}\| \coloneqq \infty$ for the operator norm.

\begin{theorem}\label{thm:applicabilityConditions}
Let $ p \in [1,\infty]$, $\alpha \in [0,1]$ be irrational, $\theta \in [0,1)$, $\lambda\in\RR$,
and $H_{\lambda, \alpha, \theta}$ be the aperiodic Schrödinger operator given by~\eqref{eq:SchrodingerOp}. For any $k \in \ZZ$, we set $\theta_k\coloneqq k\alpha\bmod 1$.
Then the following are equivalent: 
 \begin{enumerate}[label=(\roman*)]
   \item\label{it:fsmH}  The finite section method for $H_{\lambda, \alpha, \theta}$ is applicable.
   \item\label{it:Hpsi}  The operators $(H_{\lambda, \alpha, \psi})_+$ and $(\widetilde{H}_{\lambda, \alpha, \psi})_+$  are invertible for every $\psi\in [0,1)$.
     \item\label{it:HpsiUnif}  The operators $(H_{\lambda, \alpha, \theta_k})_+$ are invertible for all $k\in\ZZ$ and their inverses are uniformly bounded.
     \item\label{it:Htilde+Unif} The operators $(\widetilde{H}_{\lambda, \alpha, \theta_k})_+$ are invertible for all $k\in\ZZ$ and their inverses are uniformly bounded.
     \item\label{it:HStable} 
     For all $k\in\ZZ$, the sequence $((H_{\lambda, \alpha_m, \theta_k})_+)_{m\in\NN}$ is stable and
 \begin{equation}\label{eq:uniformStability}
  \limsup_{m\rightarrow \infty} \big( \max_{k \in \ZZ} \big\|(H_{\lambda, \alpha_m, \theta_k})_+^{-1}\big\| \big)
 < \infty.
 \end{equation} 
 \end{enumerate}
\end{theorem}
\begin{proof}
  The equivalence~\ref{it:fsmH}
  $\Leftrightarrow$\ref{it:Hpsi} is exactly Proposition~\ref{prop:FSMiffPlus}. 
  From Proposition~\ref{prop:onlyOrbit} and the subsequent Remark~\ref{rem:onlyOrbit2}, we know that the invertibility of all  $(H_{\lambda, \alpha, \psi})_+$ and $(\widetilde{H}_{\lambda, \alpha, \psi})_+$ for $\psi\in[0,1)$ can be checked by only looking at $\theta_k$ for $k\in\ZZ$ (paying the price of having the uniform boundedness condition). 
  Therefore, we have~\ref{it:Hpsi}$\Leftrightarrow$\ref{it:HpsiUnif} $\Leftrightarrow$\ref{it:Htilde+Unif}. 
  Recall that for the operators $(H_{\lambda, \alpha, \theta_k})_+$ and $(\widetilde{H}_{\lambda, \alpha, \theta_k})_+$, $k\in\ZZ$ we have the approximation result Proposition~\ref{prop:specApproxPlus}. 
  It gives that~\ref{it:HpsiUnif} and~\ref{it:Htilde+Unif} together are equivalent to
\begin{equation}\label{eq:uniformStability2}
 \sup_{k\in\ZZ} \big(\limsup_{m\rightarrow \infty}\big\|(H_{\lambda, \alpha_m, \theta_k})_+^{-1}\big\| \big)< \infty.
 \end{equation}
Therefore it only remains to show that
  \begin{equation*}
	 \sup_{k\in\ZZ} \big(\limsup_{m\rightarrow \infty}\big\|(H_{\lambda, \alpha_m, \theta_k})_+^{-1}\big\| \big)= \limsup_{m\rightarrow \infty} \big( \max_{k \in \ZZ} \big\|(H_{\lambda, \alpha_m, \theta_k})_+^{-1}\big\| \big)	\, .
  \end{equation*}
In order to prove this, note that $H_{\lambda, \alpha_m, \theta_k}$ is a $q_m$-periodic Schrödinger operator
and hence we find that 
  \begin{align*}
  \sup_{k \in \ZZ} \big\|(H_{\lambda, \alpha_m, \theta_k})_+^{-1}\big\|=
  \max_{k \in \ZZ} \big\|(H_{\lambda, \alpha_m, \theta_k})_+^{-1}\big\|\,.
  \end{align*}
 Moreover, observe that
  \begin{equation*}
  \sup_{k\in\ZZ} \big(\limsup_{m\rightarrow \infty}\big\|({H}_{\lambda, \alpha_m, \theta_k})_+^{-1}\big\| \big) \le  \limsup_{m\rightarrow \infty} \big( \sup_{k \in \ZZ} \big\|(H_{\lambda, \alpha_m, \theta_k})_+^{-1}\big\| \big)\,.
  \end{equation*} 
  
  For the other inequality, we use that the inverse of the lower norm is the norm of the inverse operator. 
  Therefore, we want to show
  \begin{equation}\label{eq:liminf_swap}
   \inf_{k\in\ZZ} \big(\liminf_{m\rightarrow\infty} \nu(H_{\lambda, \alpha_m, \theta_k})_+\big) \le \liminf_{m\rightarrow\infty}( \inf_{k\in\ZZ}\nu(H_{\lambda, \alpha_m, \theta_k})_+) \,.
  \end{equation}

Let $\delta$ denote the right-hand side of \eqref{eq:liminf_swap}.  
Then for all $N\in\NN$ there are integers $k(N)$ and $m(N)$ such that $q_{m(N)-1}>k(N)$ and
  \begin{equation*}
  \nu\big((H_{\lambda, \alpha_{m(N)}, \theta_{k(N)}})_+\big) \le \delta +  \frac{1}{N}\,.
  \end{equation*}
  Then Lemma~\ref{lem:potential_aperiodic_theta} yields $v_{\alpha_{m(N)},\theta_{k(N)}}(l) = v_{\alpha, \theta_{k(N)}}(l)$ for $-q_{m(N)}+1 \le l \le q_{m(N)}-2$.
  
Arguing as in the proof of Proposition~\ref{prop:onlyOrbit}, we obtain a $B\in\Lla$ with 
  \begin{equation*}
\delta+ \frac{1}{N_j}\ge  \nu(\big(H_{\lambda, \alpha_{m(N_j)}, \theta_{k(N_j)}})_+\big) \xrightarrow{j\rightarrow \infty}
  \nu(B_+) )  \, .
  \end{equation*}
  Hence, $\nu(B_+)  \le \delta$.
Now  Proposition~\ref{prop:specApproxPlus} and Proposition~\ref{prop:onlyOrbit} give
\begin{equation}\label{eq:limInfDelta}
 \inf_{k \in \ZZ}  \big( \liminf_{m \rightarrow \infty} \nu((H_{\lambda, \alpha_m, \theta_k})_+) \big) \le \delta\,,\end{equation}
which proves \eqref{eq:liminf_swap}.
\end{proof}

\begin{remark}\label{rem:applicabilityTilde}
The applicability of the FSM for $\widetilde{H}_{\lambda, \alpha, \theta}$ is equivalent to the applicability of the FSM for $H_{\lambda, \alpha, \theta}$, and therefore to any of the conditions~\ref{it:fsmH}--\ref{it:HStable}.
Furthermore, in part~\ref{it:HStable}, the operators $(H_{\lambda, \alpha, \theta_k})_+$ can be replaced by $(\widetilde{H}_{\lambda, \alpha, \theta_k})_+$.
\end{remark}

\begin{remark}\label{rem:FSMNorms}
Let $A_n$ denote the finite submatrices that are obtained by applying the FSM to $H_{\lambda, \alpha, \theta}$. 
Combining  \cite[Theorem~6.6]{Lindner.2016} with the above proof, one obtains 
  \begin{align*}
\|A_n^{-1}\| \le  \sup_{k\in\ZZ} \big(\limsup_{m\rightarrow \infty}\big\|(H_{\lambda, \alpha_m, \theta_k})_+^{-1}\big\| \big)
  \end{align*}
if the FSM is applicable.
\end{remark}

  \begin{corollary}\label{cor:fsm-ap}
    Let $p\in[1,\infty]$, $\lambda\in\RR$, $\alpha\in[0,1]$ be irrational with rational approximants
    $\alpha_m = \frac{p_m}{q_m}$, $\theta \in [0,1)$,
    and $H_{\lambda, \alpha, \theta}$ be the aperiodic Schrödinger operator given by~\eqref{eq:SchrodingerOp}.
    Then, the finite section method for~$H_{\lambda, \alpha, \theta}$ is applicable if and only if
    there exists~$m_0 \geq 0$ such that
    \begin{equation}\label{eq:traceCondition2} 
      \left|\tr\big(M_{\lambda,\alpha,0}(q_{m_0}) \big)\right| 
      > 2 
      \quad \text{and} \quad 
      \left| \tr \big( M_{\lambda,\alpha,0}(q_{m_0+1}) \big)\right| 
      > 2
    \end{equation}
    and
    \begin{equation}\label{eq:spectralDistance}
      \liminf_{m \rightarrow \infty} \dist(0,G_m) \neq 0
    \end{equation}
    with
    \begin{equation*}
      G_m \coloneqq
      \bigcup_{k=1}^{q_m} 
      \Big\{
        E \in \sigma\Big( \big( (A_m)_{ij}\big)_{i,j = k}^{q_m + k - 2}\Big) 
        : \left|\det\Big( \big((A_m)_{ij} - E\delta_{ij}\big)_{i,j = k}^{q_m + k - 1}\Big)\right| < 1
      \Big\}
    \end{equation*}
    where $A_m$ denotes the two-sided infinite matrix representing $H_{\lambda, \alpha_m, 0}$\,.
  \end{corollary}
\begin{proof}
We start by proving the statement for $p=2$.
By Theorem~\ref{thm:applicabilityConditions} part~\ref{it:HStable}, the finite section method for $H_{\lambda, \alpha, \theta}$ is applicable if and only if
for  all $k\in\ZZ$ the sequence $((H_{\lambda, \alpha_m, \theta_k})_+)_{m\in\NN}$ is stable and~\eqref{eq:uniformStability} holds.

From~\eqref{eq:traceCondition2} it follows there is an $\varepsilon>0$ such that for all $E$ with $|E|<\varepsilon$ we have 
\begin{equation*}
 \abs{\tr(M_{\lambda, \alpha_{m_0},\theta}(q_{m_0},E))} > 2 \quad \textnormal{and} \quad \abs{\tr(M_{\lambda, \alpha_{m_0+1},\theta}(q_{m_0+1},E))} > 2 \,,
\end{equation*}
Due to Proposition~\ref{prop:trace_condition1} and Remark~\ref{rem:traces1} we have that~\eqref{eq:traceCondition2} is equivalent to
$ E\not\in \sigma(H_{\lambda, \alpha_m, 0})$ for $m\geq m_0$. 
By Proposition~\ref{prop:lowerNormsEqual} and the fact that $\nu(A)=\dist(0,\sigma(A))$ for a self-adjoint operator $A$, this is equivalent to 
  \begin{equation*}
    0\not\in \limsup_{m\in \NN} \bigcup_{k\in\ZZ}  \sigma(H_{\lambda, \alpha_m, \theta_k})\,.
  \end{equation*}
Further, by Proposition~\ref{prop:oneSidedSpecRep}, $E\in\sigma((H_{\lambda, \alpha_m, \theta_k})_+)$ if and only if 
$E\in \sigma(H_{\lambda, \alpha_m, \theta_k})$ or both
 $E \in \sigma(((H_{\lambda, \alpha_m, \theta_k})_{i,j})_{i,j = 0}^{q_m - 2})$ 
 and $ |\det(((H_{\lambda,\alpha_m,\theta_k})_{i,j}- E\delta_{i,j})_{i,j = 0}^{q_m-1})| < 1$.
Now let the set $G_m^{\prime}$ be defined by
  \begin{equation*}
   \bigcup_{k \in \NN} 
		\Big\{
		E \in \sigma(((H_{\lambda, \alpha_m, \theta_k})_{i,j})_{i,j = 0}^{q_m - 2}) : |\det(((H_{\lambda,\alpha_m,\theta_k})_{i,j} - E\delta_{i,j})_{i,j = 0}^{q_m-1})| < 1
		\Big\}\,.
  \end{equation*}
Observe that, for all $k\in\ZZ$, the operators $H_{\lambda,\alpha_m,\theta_k}$ are $q_m$-periodic Schrödinger operators over the same alphabet. 
Even more is true: their corresponding Sturmian words $v_{\alpha, \theta_k}$ only differ by a circular shift.
Hence in fact $G_m^{\prime} = G_m$.
Consequently,~\eqref{eq:traceCondition2} and~\eqref{eq:spectralDistance} together are equivalent to
\begin{equation}
	\label{eq:LimsupinProofofMain1}
0\not\in \limsup_{m\in \NN} \bigcup_{k\in\ZZ} \sigma((H_{\lambda, \alpha_m, \theta_k})_+).
\end{equation}
Now, we use that $\|A^{-1}\| =  \dist(0,\sigma(A))$ for an invertible self-adjoint operator $A$ and obtain that 
\eqref{eq:LimsupinProofofMain1}
is equivalent to~\ref{it:HStable} in Theorem~\ref{thm:applicabilityConditions}.
 
 Having shown this, the statement indeed holds for all $p\in[1,\infty]$ because, by Proposition~\ref{prop:kurbatov}, both the set $G_m$ and the applicability of the FSM are independent of $p\in[1,\infty]$. 
 \end{proof}

\begin{figure}[htbp]
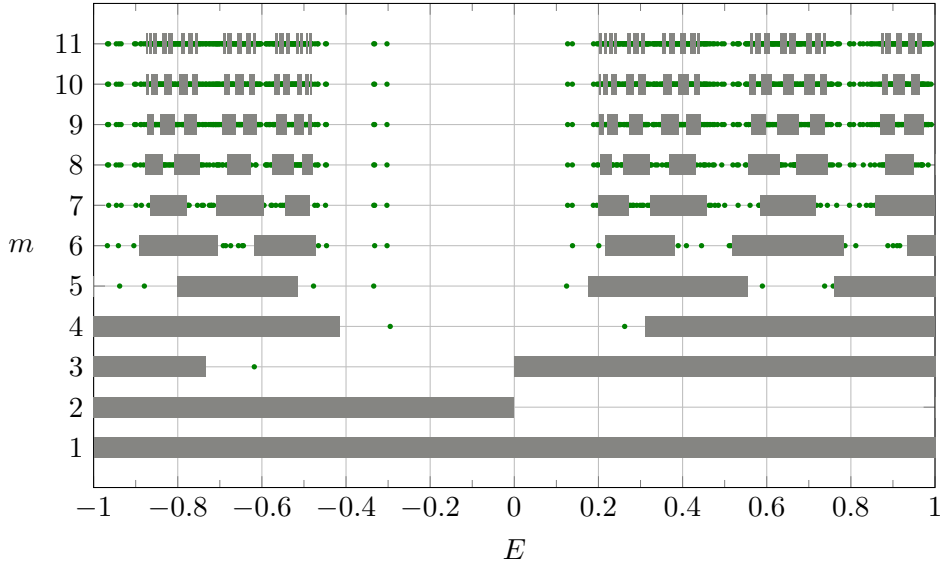
 
    \def\drawBand{\fill[bandSpecCol]}
\def\drawPoint{\fill[pointSpecCol]}

     \caption{Approximation of the spectrum of the one-sided Fibonacci Hamiltonian and its limit operators restricted to energies $E \in [-1,1]$. 
  The band spectrum of the two-sided periodic Schrödinger operators $H_{1, \alpha_m, \theta}$ is represented by the colour grey.
  The point spectrum of the one-sided periodic Schrödinger operators, given by $G_m$ in Corollary~\ref{cor:fsm-ap}, is represented by the colour green.
  }\label{fig:FibHamScep}
\end{figure}

\begin{example}\label{ex:FibHamFinal}

In Figure~\ref{fig:FibHamScep}, the spectra of the periodic approximations are displayed for the case~$\lambda=1$, $\alpha=(\sqrt{5}-1)/2$ and $\theta =0$.
As the applicability of the finite section method is concerned about excluding the spectral value $0$, the visualisation is restricted to parts of the spectrum that reside in the interval~$[-1,1]$.
The grey bars display the spectra of~$H_{\lambda, \alpha_m, 0}$, which coincide with their essential spectra due to Proposition~\ref{prop:essSpecEqualsSpecnew}.  
Furthermore, recall from Proposition~\ref{prop:lowerNormsEqual} that~$\sigma(H_{\lambda, \alpha_m, \theta}) = \sigma(H_{\lambda, \alpha_m,0})$ for all offsets~$\theta \in [0,1)$.
The green points represent the sets $G_m$ which contain the discrete spectra of the one-sided periodic operators $(H_{\lambda, \alpha_m, \theta})_+$ according to Proposition~\ref{prop:oneSidedSpecRep}.
Consequently, the whole spectrum of~$(H_{\lambda, \alpha_m, \theta})_+$ for fixed~$m$ is given by the union of the grey bands and the green points.

The figure allows to draw conclusions for the spectrum of the two-sided aperiodic case in the sense of Proposition~\ref{prop:SpurKriterium}:
If an energy~$E \in \RR$ is not contained in a grey band for two consecutive values of $m$, then it will not be contained in any grey band for larger~$m$, and also not in the spectrum of~$H_{\lambda, \alpha, \theta}$.
For instance from observing the bands for $m = 4$ and $m = 5$, one can already deduce that~$0\not\in \sigma (H_{\lambda,\alpha, \theta})$ and hence~\eqref{eq:traceCondition2} is satisfied.
For the spectrum of~$(H_{\lambda, \alpha, \theta})_+$, however, the results of Section~\ref{sec:specHltheta} do not guarantee that no green points approach $0$ for larger~$m$.
This observation will be further pursued and resolved in Section~\ref{sec:checkingApplicability}.
\end{example}

\subsection{Checking the Applicability of the FSM Numerically}\label{sec:checkingApplicability}
In order to examine the applicability of the FSM for aperiodic Schrödinger operators, we have to check the invertibility of all one-sided Schrödinger operators with the same slope $\alpha$, see Theorem~\ref{thm:applicabilityConditions} above.
In this section, we describe a method to verify the invertibility of all these operators, see Algorithm~\ref{alg:lowerNormsl2} and the code in~\cite{Gabel.2021} for an implementation.  
Again, we investigate the invertibility via the lower norm. 
Recall from~\eqref{eq:lowerNormPlus} that the approximate lower norm of a one-sided Schrödinger operator~$A_+$ is attained either with an element supported at~$0$ or it is also attained by the two-sided operator~$A$. 
We will see that, for the applicability of the FSM, we only need to look at elements with support at~$0$.
Hence, for $D\in\NN$, we write $[A]_{D}\colon \CC^{D+1}\rightarrow \CC^{D+2}$ for the restriction $A_{+}\chi_{[0,D]}$ with the natural embeddings.
As in Section~\ref{sec:specHltheta}, we set $\theta_k\coloneqq k\alpha \bmod 1$ for $k\in\ZZ$. 

As the following results will show, a key tool to investigate the applicability of the FSM for $H_{\lambda, \alpha, \theta}$ is the set
 $\Gamma_D \coloneqq \{ [H_{\lambda,\alpha,\theta_k}]_{D} : k \in \ZZ\} \subset \CC^{(D+2) \times (D+1)}$.
Since $v_{\alpha, \theta}$ has only $D+2$ different subwords of length $D+1$, see~\cite{Berstel.1996}, we find that $|\Gamma| = D + 2$.
See~\cite{Perrin.2012} for an efficient algorithm to generate these subwords.

The foundation of our methods is a quantified version of Lemma~\ref{lem:lowerNorm}. 
While explicit bounds for the truncation error can be found for all $p\in[1, \infty]$ in~\cite{Lindner.2014},
 the sharpest results are known for the case $p=2$, 
 see~\cite{Chandler-WildeLindner.2023.1,Chandler-WildeLindner.2023.2, Chonchaiya.2010}.
 Thus, in all what follows, we set $p=2$.
\begin{lemma}\label{lem:lowerNormOneSided}
	Let $\alpha\in[0,1]$ be irrational, $\theta\in[0,1)$, and $\lambda\in\RR$. 
  For $D\in\NN$ and $\varepsilon_D>0$  satisfying 
	\begin{equation}\label{eq:lowerNormIneq2}
		\varepsilon_D > 4 \sin\frac \pi {2D+2},
	\end{equation} 
   we have that
	\begin{equation*}
	\inf_{k\in\ZZ}\nu((H_{\lambda,\alpha,\theta_k})_+) \ge \min_{B\in\Gamma_D}\nu(B)  -\varepsilon_D.
	\end{equation*}
\end{lemma}
\begin{proof} For $k\in\ZZ$, let us write $A_k\coloneqq H_{\lambda,\alpha,\theta_k}$. 
  Using the bounds obtained in~\cite[Corollary~4.4]{Chonchaiya.2010}, we get the lower norm estimate 
  $\nu((A_k)_+)\ge \nu_D((A_k)_+)-\varepsilon_D$.
	Now we take the infimum over $k \in \ZZ$ and apply Lemma~\ref{lem:subwordsAndLowerNorms}\ref{it:decomposition}:
        \begin{equation*}
          \inf_{k\in\ZZ}\nu((A_k)_+) \ge  \inf_{k\in\ZZ} \min  \left\{\nu([A_k]_{D}),  \nu_D(A_k)\right\}  -\varepsilon_D \,.
        \end{equation*}
	It is easy to verify that, for all $j\in\ZZ$,
        \begin{equation*}
          \nu_D(A_j)\ge \inf_{k\in\ZZ}  \nu([A_k]_{D}) = \min_{B\in\Gamma_D} \nu(B)
        \end{equation*}
	holds
\end{proof}

Let us put the previous result into the context of the FSM.

\begin{theorem}\label{thm:numerikFSM}
	Let $\alpha\in [0,1]$ be irrational, $\theta\in[0,1)$, and $\lambda\in\RR$.  
  Then the FSM is applicable to $H_{\lambda, \alpha,\theta}$ 
  if and only if there are $D\in\NN$ and $\varepsilon_D$ satisfying~\eqref{eq:lowerNormIneq2} such that 
        \begin{equation}\label{eq:numFSM_condition}
          \min_{B\in\Gamma_D} \nu(B) > \varepsilon_D \,.
        \end{equation}
	In that case, we have that 
        \begin{equation}\label{eq:FSMBound}
          \limsup_{n\rightarrow\infty} \big\| ((H_{\lambda,\alpha,\theta})_{i,j=0}^{n})^{-1}\big \| 
          \le \frac{1}{\min\limits_{B\in\Gamma_D}\nu(B)  - \varepsilon_D}\,.
        \end{equation}
\end{theorem}
\begin{proof}
	We use Theorem~\ref{thm:applicabilityConditions}, namely the equivalence \ref{it:fsmH}$\Leftrightarrow$\ref{it:HpsiUnif}. Hence, the FSM is applicable to  $H_{\lambda, \alpha,\theta}$ if and only if 
	$ \inf_{k\in\ZZ}  \nu((H_{\lambda, \alpha, \theta_k})_+) >0 $.
	Here we replaced the norms of the inverses by the lower norms.

Now we can apply Lemma~\ref{lem:lowerNormOneSided} and obtain
        \begin{equation*}
          \inf_{k\in\ZZ} \nu([H_{\lambda,\alpha,\theta_k}]_{D}  ) -\varepsilon_D
           \le \inf_{k\in\ZZ}\nu((H_{\lambda,\alpha,\theta_k})_+)    
           \le    \inf_{k\in\ZZ} \nu([H_{\lambda,\alpha,\theta_k}]_{D}  )\, . 
        \end{equation*}	
        On the one hand, this shows that the FSM is applicable if  $\nu([H_{\lambda,\alpha,\theta_k}]_{D} ) >\varepsilon_D $. 
        On the other hand, if the FSM is applicable, then $ \inf_{k\in\ZZ} \nu([H_{\lambda,\alpha,\theta_k}]_{D}  )$ is bounded from below by $\inf_{k\in\ZZ}\nu((H_{\lambda,\alpha,\theta_k})_+) >0$. 
Then, for sufficiently large $D$, the corresponding $\varepsilon_D$ becomes smaller than this bound. 
Hence, $\nu([H_{\lambda,\alpha,\theta_k}]_{D} ) >\varepsilon_D$. 
This concludes the proof of the first statement.

To prove~\eqref{eq:FSMBound}, it can be found in~\cite[Theorem~6.6]{Lindner.2016} that the bound on the inverses of the finite sections for an operator $A$ is equal to the supremum of the norms of the inverses of the following operators:

\begin{inparaenum}[(a)]
  \item\label{it:opA} $\,A\,,$ \qquad
\item\label{it:opB-} $\,B_-$ for all $B\in  \Lim_+(A) \,,$ \qquad
\item\label{it:opC+} $\,C_+$ for all $C\in  \Lim_-(A) \,.$
\end{inparaenum}
Again, we use that the norms of the inverses are the inverses of the lower norms.	
The identity \eqref{eq:FSMBound} then follows from above considerations.
\end{proof}

As a consequence of the above, Algorithm~\ref{alg:lowerNormsl2} can be used to check the applicability of the FSM.
\begin{algorithm}
	\caption{
    Algorithm to check the applicability of the FSM for $H_{\lambda, \alpha, \theta}$ based on Theorem~\ref{thm:numerikFSM}.
    For an efficient way of implementing \texttt{generateAllWords}, see~\cite{Perrin.2012}. 
    The function \texttt{constructMatrix} returns the  $(D+2) \times (D+1)$ matrix with $w$ on its main diagonal and $1$ on the first sub- and superdiagonal 
    and \texttt{minSVD} returns the smallest singular value.    
  }\label{alg:lowerNormsl2}
	\begin{algorithmic}[1]
		\Procedure{checkApplicability}{$\lambda$, $\alpha$, $D$}
		\State $\varepsilon_D = 4 \sin\frac \pi {2(D+2)}$

		\State words = generateAllWords($\lambda$, $\alpha$, $D$)
		\Comment generates all subwords of length $D$
		\State $N=\emptyset$
		\For{$w\in$ words}
		\State$H$ = constructMatrix$(w)$
		\State $\nu$ = minSVD($H$)
		\State$N$ = $N \cup$ $\{\nu \}$
		
		\EndFor
		
		\If{$\min(N) \le \varepsilon_D$}
		\State	\Return false	
		\Else \State
		\Return true
		\EndIf
		\EndProcedure
	\end{algorithmic}
\end{algorithm}

In the following example, we use Algorithm~\ref{alg:lowerNormsl2} to recover and extend the results of~\cite{Lindner.2018} where a different approach was used. 
In particular, we show that the FSM for the Fibonacci Hamiltonian is applicable.

\begin{example} \label{ex.FiboHam_num}
	We consider the case $\lambda=1$ and set $D=100$. 
        With $\varepsilon_D \coloneqq 4\sin \frac\pi{202} \approx 0.062$, we check whether 
        $\nu((H_{\lambda,\alpha,\theta_k})_{i,j=0}^{D})= \big\|\big((H_{\lambda,\alpha,\theta_k})_{i,j=0}^{D}\big)^{-1}\big\|^{-1}> \varepsilon_D$
        holds for every $k\in\ZZ$. Take any $\alpha$ with a continued fraction expansion
	\begin{equation}\label{eq:exampleFractions}
    \alpha = [a_1, a_2, \dots]  \quad \text{with} \quad a_1 = a_2 = \dots = a_{11} = 1.
	\end{equation}
	Using these digits in the recursion formula~\eqref{eq:recursiveWord}, we obtain a word $s_{11}$ of length $144$. 
  Since we are only interested in subwords of length $D+1=101$,
	all other $a_j$ for $j>11$ yield the same result. 
	Computing the smallest singular values, i.e.~the lower norms for all $102$ subwords of $v_{\alpha, 0}$ of length $101$, we find that
        \begin{equation*}\label{eq:LN_results}
	\min_{B\in\Gamma_D}\nu(B) \ge 0.126 > \varepsilon_D \,.
        \end{equation*}
	Therefore,~\eqref{eq:numFSM_condition} holds and we can use Lemma~\ref{lem:lowerNormOneSided} to conclude that 
        \begin{equation*}
	\inf_{k\in\ZZ}\nu((H_{\lambda,\alpha,\theta_k})_+) \ge 0.126 - \varepsilon_D \ge 0.063 >0.
        \end{equation*}
	
	By Theorem~\ref{thm:applicabilityConditions} we find that $(H_{\lambda,\alpha,\theta})_+$ is invertible for every $\theta\in [0,1)$.
	Now Theorem~\ref{thm:numerikFSM} implies that the FSM is applicable to $H_{\lambda, \alpha,\theta}$ for all $\alpha$ of the form~\eqref{eq:exampleFractions} 
  and all $\theta\in [0,1)$. 
  Moreover, due to~\eqref{eq:limsupAm}, 
  \begin{equation*}
    \limsup_{n\rightarrow\infty} \big\| ((H_{\lambda,\alpha,\theta})_{i,j=0}^{n})^{-1}\big \| \le \frac{1}{0.063} \approx 15.87\,.
  \end{equation*}

  The above computations also reveal that there is a $k\in \ZZ$ (in fact, this is the case for most $k$) such that $\nu([H_{\lambda,\alpha,\theta_k}]_D)\approx 0.127$.
  However, is clear from Proposition~\ref{prop:SpurKriterium} and Figure~\ref{fig:FibHamScep} that $\text{dist}(0, \sigma_\ess((H_{\lambda, \alpha, \theta_k}))_+) \ge 0.2$.
  Since
  \begin{equation}
    \nu((H_{\lambda, \alpha, \theta_k})_+) \le \nu([H_{\lambda,\alpha,\theta_k}]_D)  <  \text{dist}(0, \sigma_\ess((H_{\lambda, \alpha, \theta_k})_+)) ,
  \end{equation}
  we find that $(H_{\lambda, \alpha, \theta_k})_+$ has to have an eigenvalue $E_0$ in the gap around $0$ with $0.127- \varepsilon_D <|E_0|  <0.127$.
  Note how Figure~\ref{fig:FibHamScep} also shows positive eigenvalues in this range for the periodic approximants.

Figure~\ref{fig:ln_vs_eps} compares $\min_{B\in\Gamma_D}\nu(B)$ and $\varepsilon_D$ for $D=1, \dots, 100$. 
We can see that the crucial condition~\eqref{eq:numFSM_condition}, i.e.~$\min_{B\in\Gamma_D}\nu(B) > \varepsilon_D$ holds for $D \ge 49$.

  \begin{figure}[htbp] 
%
%
\definecolor{mycolor2}{rgb}{0.4,0.1,0.1}%
\definecolor{mycolor1}{rgb}{0.4,0.6,0.95}%
\begin{tikzpicture}

\begin{axis}[%
width=4.200in,
height=3.462in,
at={(0.8in,0.508in)},
scale only axis,
xmin=1,
xmax=100,
xlabel style={font=\color{white!15!black}},
xlabel={$D$},
ymin=0,
ymax=0.5,
xmajorgrids,
ymajorgrids,
axis background/.style={fill=white},
legend style={legend cell align=left, align=left, draw=white!15!black}
]
\addplot [scatter, only marks,
scatter/classes={0={mark=*,mycolor1, scale=0.9}}]
  table[row sep=crcr]{%
1	0.662153446861956\\
2	0.569972919912301\\
3	0.419462686793632\\
4	0.348492380348655\\
5	0.313845282945127\\
6	0.230939147833251\\
7	0.227845781848129\\
8	0.194376295740341\\
9	0.180822362683835\\
10	0.173711478357719\\
11	0.158639012391273\\
12	0.154473723619003\\
13	0.14947339009459\\
14	0.140754415918406\\
15	0.139978106440754\\
16	0.135125709698912\\
17	0.13341747345449\\
18	0.132229000323875\\
19	0.129990717641872\\
20	0.129846948984968\\
21	0.128806671882781\\
22	0.128443612186942\\
23	0.128209042907455\\
24	0.127768293318544\\
25	0.127659440064867\\
26	0.127506294433817\\
27	0.127249487019843\\
28	0.127223310844897\\
29	0.127069745629111\\
30	0.127015848269569\\
31	0.126976314146565\\
32	0.126901832620757\\
33	0.126883144524744\\
34	0.126855274902464\\
35	0.126807337320221\\
36	0.126802215850462\\
37	0.126772052878881\\
38	0.126761283648577\\
39	0.126753259525563\\
40	0.126737993617692\\
41	0.126736953810816\\
42	0.126729539751553\\
43	0.126726934865414\\
44	0.126725219565366\\
45	0.126721989642281\\
46	0.126721190572889\\
47	0.126720053580069\\
48	0.126718139889095\\
49	0.126717942765022\\
50	0.126716787015491\\
51	0.12671638033516\\
52	0.126716080783212\\
53	0.126715515618454\\
54	0.126715477469656\\
55	0.126715205989875\\
56	0.126715110976194\\
57	0.126715048554906\\
58	0.126714931238318\\
59	0.12671490226237\\
60	0.126714861073904\\
61	0.12671479183863\\
62	0.12671478471402\\
63	0.126714742964849\\
64	0.126714728284265\\
65	0.126714717474378\\
66	0.126714697086995\\
67	0.126714691967311\\
68	0.126714684313842\\
69	0.126714671132188\\
70	0.126714669720575\\
71	0.126714661404509\\
72	0.126714658432561\\
73	0.126714656216098\\
74	0.126714651996759\\
75	0.12671465170913\\
76	0.126714649658097\\
77	0.12671464893727\\
78	0.126714648462492\\
79	0.126714647568357\\
80	0.126714647347126\\
81	0.126714647032296\\
82	0.126714646502344\\
83	0.126714646447749\\
84	0.126714646127638\\
85	0.126714646014992\\
86	0.126714645932016\\
87	0.126714645775458\\
88	0.12671464573613\\
89	0.126714645677327\\
90	0.126714645576022\\
91	0.126714645565171\\
92	0.126714645501238\\
93	0.126714645478386\\
94	0.126714645461342\\
95	0.126714645428894\\
96	0.126714645426682\\
97	0.126714645410907\\
98	0.126714645405363\\
99	0.126714645401711\\
100	0.126714645394834\\
};
\addlegendentry{$\min_{\Gamma_D}\nu(B)$}

\addplot [scatter, only marks,
scatter/classes={0={mark=10-pointed star,mycolor2, scale=1}}]
  table[row sep=crcr]{%
1	2.82842712474619\\
2	2\\
3	1.53073372946036\\
4	1.23606797749979\\
5	1.03527618041008\\
6	0.890083735825258\\
7	0.780361288064513\\
8	0.694592710667721\\
9	0.625737860160923\\
10	0.569259353093141\\
11	0.522104768880206\\
12	0.482146721021292\\
13	0.447857904413231\\
14	0.418113853070614\\
15	0.392068561318242\\
16	0.369073437853208\\
17	0.348622970990633\\
18	0.330317381889329\\
19	0.31383638291138\\
20	0.298920374345697\\
21	0.285356732796929\\
22	0.272969653458684\\
23	0.261612516920572\\
24	0.251162078117253\\
25	0.241513989689144\\
26	0.232579315641903\\
27	0.224281788948767\\
28	0.21655563434167\\
29	0.209343824971775\\
30	0.202596675354851\\
31	0.196270697309672\\
32	0.190327663294969\\
33	0.184733834582958\\
34	0.17945932140206\\
35	0.174477549461344\\
36	0.169764812784593\\
37	0.165299896995253\\
38	0.161063760437661\\
39	0.157039263036274\\
40	0.153210934760141\\
41	0.149564777105302\\
42	0.146088092230635\\
43	0.142769335355922\\
44	0.139597986810004\\
45	0.136564440743872\\
46	0.133659908030698\\
47	0.130876331287105\\
48	0.128206310286621\\
49	0.125643036312513\\
50	0.123180234224681\\
51	0.120812111203555\\
52	0.118533311290239\\
53	0.116338874972446\\
54	0.114224203174785\\
55	0.112185025103476\\
56	0.110217369472648\\
57	0.108317538704538\\
58	0.106482085751099\\
59	0.104707793231493\\
60	0.102991654619954\\
61	0.101330857252752\\
62	0.0997227669522915\\
63	0.0981649140916492\\
64	0.0966549809445292\\
65	0.0951907901844301\\
66	0.0937702944130418\\
67	0.0923915666120044\\
68	0.0910527914244227\\
69	0.0897522571832197\\
70	0.0884883486127476\\
71	0.0872595401382445\\
72	0.0860643897448888\\
73	0.0849015333344941\\
74	0.0837696795334278\\
75	0.0826676049102172\\
76	0.0815941495656214\\
77	0.0805482130617613\\
78	0.0795287506602807\\
79	0.0785347698425132\\
80	0.0775653270872975\\
81	0.0766195248844564\\
82	0.075696508964076\\
83	0.0747954657236091\\
84	0.0739156198365196\\
85	0.0730562320276992\\
86	0.0722165970022418\\
87	0.0713960415153842\\
88	0.0705939225725144\\
89	0.069809625749134\\
90	0.069042563621554\\
91	0.0682921742998989\\
92	0.0675579200557244\\
93	0.0668392860372038\\
94	0.0661357790654338\\
95	0.0654469265059471\\
96	0.064772275210006\\
97	0.064111390520695\\
98	0.0634638553392317\\
99	0.0628292692472827\\
100	0.0622072476814035\\
};
\addlegendentry{$\varepsilon_D$}

\end{axis}
\end{tikzpicture}%
         \caption{Comparison of the quantities involved in condition~\eqref{eq:numFSM_condition} 
        with $\alpha$ and $\lambda$ from Example~\ref{ex:FibHamFinal} and $D\in \{1, \dots, 100\}$.
        It can be seen that the conditions of Theorem~\ref{thm:numerikFSM} hold for $D \ge 49$.
      }\label{fig:ln_vs_eps}
  \end{figure}
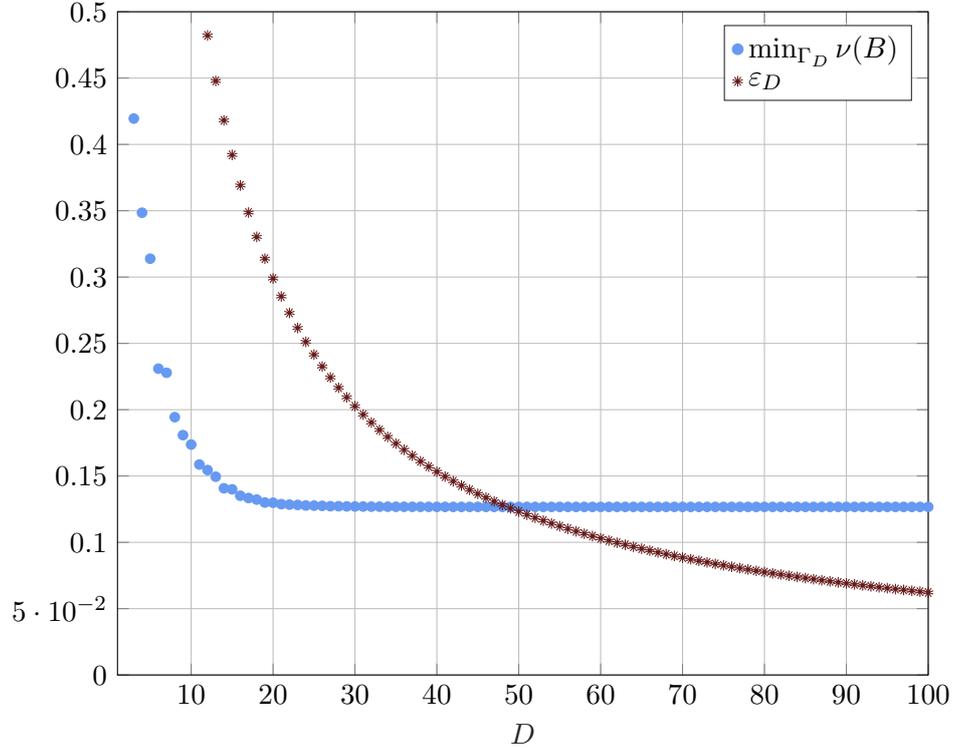

\end{example}

\section*{Acknowledgement}
The authors would like to thank Alexander Haupt for providing stimulating discussions about Sturmian words.

\typeout{get arXiv to do 4 passes: Label(s) may have changed. Rerun}

\end{document}